\documentclass[12pt]{amsart}
\usepackage{amsmath,amssymb,amscd,array, mathrsfs }

\usepackage{amsmath,amscd,amssymb,amsthm,array}
\usepackage{color}

\usepackage{amsmath,amssymb,mathrsfs,amsthm, tikz-cd,mathrsfs}
\usepackage{comment}
\def\url#1{\expandafter\s

\tring\csname #1\endcsname}

\def\mmat #1,#2,#3,#4,{\text{\small\arraycolsep=3pt $
\begin{pmatrix}#1&#2\\#3&#4\end{pmatrix}$}}

\usepackage{hyperref}
\usepackage{capt-of}

\usepackage{multirow}
\usepackage[all]{xy}

\usepackage{comments}
\newComments\SBe{Said}{blue}
\newComments\SBo{Sofiane}{blue}
\newComments\AM{Nacer}{blue}
\newComments\DL{DL}{red}
\newComments\QEh{QEh}{blue}

\def\mmat #1,#2,#3,#4,{\text{\small\arraycolsep=3pt $
\begin{pmatrix}#1&#2\\#3&#4\end{pmatrix}$}}

\usepackage{lscape}
\usepackage{tikz-cd}
\usepackage{enumerate}

\usepackage{DLdef1}
\usepackage{multicol}

\def\mmat #1,#2,#3,#4,{\text{\small\arraycolsep=3pt $
\begin{pmatrix}#1&#2\\#3&#4\end{pmatrix}$}}

\renewcommand {\ssbegin}[2][*]
 {\refstepcounter{subsection}%
\if#1*
\addcontentsline{toc}{subsection}{\thesubsection.\hskip 1pc #2}%
\else
\addcontentsline{toc}{subsection}{\thesubsection.\hskip 1pc #2. #1}%
\fi
 \def \secno {\gdef \secno {}{\ssecfont
\thesubsection.\hskip 2ex}%
 }%
 \begin{#2}}

\renewcommand {\sssbegin}[2][*]
 {\refstepcounter{subsubsection}
\if#1*
\addcontentsline{toc}{subsubsection}{\thesubsubsection.\hskip 1pc #2}%
\else
\addcontentsline{toc}{subsubsection}{\thesubsubsection.\hskip 1pc #2. #1}
\fi
 \def \secno {\gdef \secno {}{\ssecfont \thesubsubsection.\hskip 2ex}%
 }%
 \begin{#2}}

\renewcommand {\parbegin}[2][*]
 {\refstepcounter{paragraph}
\if#1*
\addcontentsline{toc}{paragraph}{\theparagraph.\hskip 1pc #2}%
\else
\addcontentsline{toc}{paragraph}{\theparagraph.\hskip 1pc #2. #1}
\fi
 \def \secno {\gdef \secno {}{\ssecfont \theparagraph.\hskip 2ex}%
 }%
 \begin{#2}}

\newcommand {\sort}{{\rm{sort}}}

\rmnameii{etr}{etr}
\rmnameii{evv}{ev}

\DeclareMathOperator{\K}{\mathbb{K}}
\newcommand{\Z}{\mathbb{Z}}

\newcommand{\black}{\color{black}}

\newcommand{\inc}{\mathrm{inc}}

\begin{document}

\title[Restricted Lie-Rinehart superalgebras]{The Superization of Hochschild's Lemma and Restricted Lie-Rinehart Superalgebras}

\author[S. Bouarroudj]{Sofiane Bouarroudj}

\address{Division of Science and Mathematics, New York University Abu Dhabi, P.O. Box 129188, Abu Dhabi, United Arab Emirates.}
\email{sofiane.bouarroudj@nyu.edu}

\author[Q. Ehret]{Quentin Ehret}
\address {Division of Science and Mathematics, New York University Abu Dhabi, P.O. Box 129188, Abu Dhabi, United Arab Emirates.}
\email{qe209@nyu.edu}

\author[A. Makhlouf]{Abdenacer Makhlouf}
\address{Université de Haute-Alsace, IRIMAS UR 7499, F-68100 Mulhouse, France.}
\email{abdenacer.makhlouf@uha.fr}

\author[N. Shyntas]{Nurtas Shyntas}
\address {Division of Science and Mathematics, New York University Abu Dhabi, P.O. Box 129188, Abu Dhabi, United Arab Emirates.}
\email{ns5493@nyu.edu}

\thanks{S.B. and Q.E. were supported by the grant NYUAD-065}



\begin{abstract}
The main goal of this paper is to introduce the notion of  restricted Lie-Rinehart superalgebra over a field of characteristic $p>2$,  motivated by a generalization of Hochschild’s lemma to the super setting. We extend Schauenburg’s proof of Hochschild’s lemma to Lie-Rinehart superalgebras and we prove a superized version that serves as a foundation for our construction. Building upon this, we define restricted Lie-Rinehart superalgebras, investigate their representations, construct  the semi-direct product with a restricted module, and provide several examples. Finally, we construct the corresponding universal enveloping algebra and show that this algebra satisfies the expected universal property. 
\end{abstract}


\maketitle

\thispagestyle{empty}
\setcounter{tocdepth}{2}
\tableofcontents

\section{Introduction}

\subsection{Restricted Lie algebras}

Over a field of characteristic $p>0$, the Lie algebra of derivations of an associative algebra has a richer structure than it does in characteristic zero, since the Frobenius map $d \mapsto d^p$ is an endomorphism of the space of derivations (see \cite{Ja}). This was the main motivati on for defining \textit{restricted} Lie algebras, introduced by Jacobson as a way to axiomatize this observation (see \cite{J}). Roughly speaking, a Lie algebra $L$ over a field of characteristic $p>0$ is called \textit{restricted} if it is equipped with an additional unary operation, called a $p$-map and denoted by 
$(-)^{[p]} : L \to L$ which mimics the $p$-th power in an associative algebra (see Section~\ref{def:restrictedLie}). 
Indeed, the canonical examples are constructed from associative algebras by taking the commutator bracket as the Lie product and the usual 
$p$-th power as the $p$-map. Geometrically,  Lie algebras associated with algebraic groups are restricted in positive characteristic and many properties of the algebraic group are encoded in the restricted structure of its Lie algebra, see \cite{S,SF}. Because many classical tools, such as the Killing form and Lie’s theorem, fail in positive characteristic, restricted Lie algebras in characteristic $p$ are used to recapture some of the behavior of Lie algebras over fields of characteristic zero, but they also exhibit new and different phenomena. Kostrikin and Šafarevič conjectured that every \textit{simple} restricted Lie algebra over an algebraically closed field of characteristic $p>7$ is of classical or Cartan type (see \cite{KS}). The conjecture was proved for $p>7$ by Block and Wilson (see \cite{BW}), while the case $p>3$ was settled by Premet and Strade (see \cite{PS} and references therein). The classification of restricted Lie algebras over a field of characteristic 2 and 3 remains an open problem, see \cite{S} for the full story. There is a special restricted cohomology theory for restricted Lie algebras introduced by Hochschild in \cite{H}. This cohomology was further studied in \cite{EF, EM2,EFa}, to cite a few. 

To the best of our knowledge, the superization is due to Mikhalev (\cite{Mi}). Several authors have since studied restricted Lie superalgebras in various contexts. The study of restricted simple Lie superalgebras in characteristic 3 can be found in \cite{BKLLS}.  Double extensions of restricted Lie superalgebras have been studied by \cite{BBH, BEM}. Additionally, identities on the restricted enveloping superalgebras have been investigated in \cite{Pe, U},  representations of the restricted Witt superalgebras and $\mathfrak{sl}(n|1)$ have been explored in \cite{SZ, Z}, and a cohomology theory pertaining to restricted Lie superalgebras has recently been developed in \cite{BE} superizing \cite{EF}. While an analogue of the Kostrikin-Šafarevič conjecture for restricted Lie superalgebras remains unclear, progress has been made with the classification of $p$-nilpotent ones in small dimensions \cite{BE}.

\subsection{Lie-Rinehart algebras} Lie-Rinehart algebras are algebraic generalizations of Lie algebroids. They were introduced in the work of Herz (\cite{He}), Palais (\cite{P}), and especially Rinehart (\cite{Ri}), who developed a formalism of differential forms for commutative algebras and studied the representations and cohomology of Lie-Rinehart algebras.
Roughly speaking, a Lie-Rinehart algebra is a triple $(A,L,\rho)$, where $A$ is an associative commutative algebra, $L$ is Lie algebra and an $A$-module, and $\rho:L\rightarrow \Der(A)$ is an $A$-linear Lie algebras morphism satisfying a compatibility condition that encodes the lack of $A$-linearity of the Lie bracket of $L$, see Eq. \eqref{eq:def-LR}. Recent works on Lie-Rinehart algebras include their relationship to Poisson algebras investigated in \cite{Hu}, where Huebschmann  constructed a Lie-Rinehart algebra from the module of K\"{a}hler differentials of a Poisson algebra, and used this construction to define Poisson cohomology and homology in terms of the corresponding structures of the associated Lie-Rinehart algebra.  Subsequent works by Huebschmann include \cite{Hu2}, which  relates Lie-Rinehart algebras with Gerstenhaber algebras and Batalin-Vilkovisky algebras, and \cite{Hu3}, which explores an application of Lie-Rinehart algebras to the quantization of a positive K\"{a}hler manifold with a Hamiltonian action by a compact Lie group. 

Over fields of positive characteristic, the notion of \textit{restricted} Lie-Rinehart algebras was introduced by Dokas in \cite{Do}, although earlier versions of the concept appear in \cite{Ho,Ru}. In particular, Hochschild proved an additional compatibility condition for the $A$-module of derivations of an associative algebra $A$, specific to the case of a field of characteristic $p>0$ (see \cite[Lemma 1]{Ho}), which is referred to as \textit{Hochschild's Lemma}. This result provided the main motivation for the definition of restricted Lie-Rinehart algebras, see Definition \ref{def:RLR}. In this case, the underlying Lie algebra $L$ is restricted, and the additional compatibility condition encodes the lack of $p$-homogeneity of the $p$-map with respect to the $A$-module structure of $L$. Dokas constructed the restricted universal enveloping algebra and provided a cohomology theory based on Quillen-Barr-Beck methods in \cite{Do}, whereas another approach adapted to formal deformations was proposed in \cite{Eh}. Moreover, Schauenburg gave an alternative construction of the restricted universal enveloping algebra of a restricted Lie–Rinehart algebra in \cite{Sc}, as well as a reformulated version of the Hochschild lemma.

Over a field of characteristic zero, Lie-Rinehart superalgebras appeared in \cite{Ch, Ro} and \cite{EM1}, where a classification in small dimensions and a theory of formal deformations were proposed.
More recently, the universal enveloping algebra of a Lie–Rinehart superalgebra was studied in \cite{La}, where the author established a new Poincaré-Birkhoff-Witt theorem for Lie–Rinehart superalgebras, leading to a similar result for Poisson superalgebras.
\subsection{The main results and the organization of the paper}
The main contribution of this paper is the definition of the notion of  {\it restricted} Lie-Rinehart {\it superalgebra} over a field of characteristic $p>2$ (see, Definition \ref{def:sRLR}), derived from the generalization of Hochschild’s lemma in the super setting, see Theorem \ref{thm:super-hochschild}. The superization of the notion of restricted Lie-Rinehart algebra is not a straightforward superization minding the sign rule, but it rather relies on the superization of Hochschild Lemma. 

The paper is organized as follows. In Section~\ref{sec:background}, we review the basic notions of Lie-Rinehart (super)algebras, restricted Lie algebras, and restricted Lie-Rinehart algebras. We also recall Hochschild's lemma for Lie-Rinehart algebras (see Proposition~\ref{prop:hochschild}). Section~\ref{sec:super-hochschild} is devoted to the superization of Hochschild's lemma, stated in Theorem~\ref{thm:super-hochschild}. Our approach extends Schauenburg's proof (see~\cite{Sc}) to the super setting and relies on tools from the theory of Hopf superalgebras. The proof is quite technical and is divided into four cases, depending on the parity of the elements involved. Some parts of the proof involving combinatorial methods are postponed to Appendix~\ref{sec:gamma2pp}. In Section~\ref{sec:sRLR}, we briefly recall the definition of restricted Lie superalgebras and introduce the notion of restricted Lie-Rinehart superalgebras in Definition~\ref{def:sRLR}, which constitutes the main contribution of this paper. Further, we introduce representations (see Definition~\ref{def:rest-rep}), and discuss the semi-direct product of a restricted Lie-Rinehart superalgebra by such a representation (see Theorem~\ref{thm:semi-direct}), and provide several examples. Finally, in Section~\ref{sec:UEA}, we define the universal enveloping algebra of a restricted Lie-Rinehart superalgebra and show that our construction satisfies the appropriate universal property (see Proposition~\ref{prop:univ-prop}).

It is presumed that when the parity of an element is written, the element is homogeneous. Throughout the text, $\K$ refers to a field of characteristic $p\geq 3$.

\section{Background}\label{sec:background}

\subsection{Lie-Rinehart superalgebras} Following \cite{Ri,Ch}, a Lie-Rinehart superalgebra is a triple $(A,L,\rho)$, where $A=A_\ev\oplus A
_\od$ is an associative supercommutative superalgebra, $(L=L_\ev\oplus L_\od,[-,-])$ is a Lie superalgebra that is also an $A$-module, and $\rho:L\rightarrow\Der(A)$ is a morphism of Lie superalgebras and of $A$-modules satisfying the \textit{Leibniz condition}:
\begin{equation}\label{eq:def-LR}
    [x,ay]=(-1)^{|x||a|}a[x,y]+\rho(x)(a)y,\quad\forall x,y\in L~\forall a\in A.
\end{equation}

Let $(A,L,\rho)$ and $(B,H,\theta)$ be Lie-Rinehart superalgebras. A morphism of Lie-Rinehart superalgebras is a pair $(\phi,\psi)$, where $\phi:A\rightarrow B$ is a morphism of associative superalgebras and $\psi:L\rightarrow H$ is a morphism of Lie superalgebras satisfying
\begin{equation}\label{eq:morph-sLR}
    \phi\circ\rho(x)(a)=\theta\bigl(\psi(x)\bigl)\bigl(\phi(a)\bigl),\quad \forall a\in A,~\forall x\in L.
\end{equation} In the case where $A=B$, one can take $\phi=\id$ and Condition \eqref{eq:morph-sLR} reduces to 
\begin{equation}
    \rho(x)(a)=\theta\bigl(\psi(x)\bigl)(a),\quad \forall a\in A,~\forall x\in L.
\end{equation}

A representation of a Lie-Rinehart superalgebra $(A,L,\rho)$ is an $A$-module $M$ together with a $A$-linear morphism of Lie superalgebras $\phi:L\rightarrow \End(M)$ satisfying
\begin{equation}\label{eq:RLRmodule}
    \phi(x)(am)=(-1)^{|x||a|}a\phi(x)(m)+\rho(x)(a)m,~\forall x\in L,~\forall a\in A,~\forall m\in M.
\end{equation} Such a pair $(\phi,M)$ is also called \textit{Lie-Rinehart module}.

In the text, we call \textit{Lie-Rinehart algebra} a Lie-Rinehart superalgebra concentrated in even degree, that is, $A_\od=L_\od=0$.
\subsection{Restricted Lie algebras}\label{def:restrictedLie} Following \cite{J}, a finite-dimensional Lie algebra $(L,[-,-])$ is called \textit{restricted} if there is a map $(-)^{[p]}:L\rightarrow L$ (called \textit{$p$-map} or \textit{$p$-structure}) satisfying
\begin{align}\label{eq:def-1}
(\lambda x)^{[p]}&=\lambda^p x^{[p]} \text{ for all } x\in L \text{ and for all } \lambda \in \K;\\\label{eq:def-2}
\ad_{x^{[p]}}&=(\ad_x)^p \text{ for all } x\in L;\\
(x+y)^{[p]}&=x^{[p]}+y^{[p]}+\displaystyle\sum_{1\leq i\leq p-1}s_i(x,y), \text{ for all } x,y\in L;\label{eq:si}
\end{align}
where the coefficients $s_i(x,y)$ can be obtained from the expansion
\begin{equation}\label{eq:def-si}
(\ad_{\lambda x+y})^{p-1}(x)=\displaystyle\sum_{1\leq i \leq p-1} is_i(x,y) \lambda^{i-1}.
\end{equation}
An explicit expression of  $\displaystyle\sum_{1\leq i\leq p-1}s_i(x,y)$ in terms of nested brackets is given by
\begin{align}\label{eq:exp-si}
\sum_{i=1}^{p-1}s_i(x,y)&=\sum_{\underset{x_{p-1}=y,~x_p=x}{x_k\in\{x,y\}}}\frac{1}{\sharp(x)}[x_1,[x_2,[...,[x_{p-1},x_p]...]]],
\end{align} where $\sharp(x)=\card\{i,~x_i=x\}\in\{1,\cdots,p-1\}$.

To investigate $p$-structures on a Lie algebra, it is useful to use the following theorem, due to Jacobson.
\sssbegin{Theorem}[\cite{J}] \label{Jac}
Let $(e_j)_{j\in J}$ be a~basis of $L$ such that there are $f_j\in L$ satisfying $(\ad_{e_j})^p=\ad_{f_j}$. Then, there exists exactly one $p$-map  $(-)^{[p]}:L\rightarrow L$ such that 
\[
e_j^{[p]}=f_j \quad \text{ for all $j\in J$}.
\]
\end{Theorem}

\sssbegin{Example}[see \cite{SF}] Let $A$ be an associative algebra. Then, $A$ can be endowed with a restricted Lie algebra structure with the bracket $[a,b]=ab-ba$ and the $p$-map $a^{[p]}=a^p$.
\end{Example}
\sssbegin{Example}[see \cite{J}] Let $A$ be an associative algebra and let $\Der(A)$ denote its space of derivations, that is, 
\begin{equation}
    \Der(A)=\bigl\{D\in\End(A),~D(ab)=D(a)b+aD(b),\quad\forall a,b\in A  \bigl\}.
\end{equation} Then, we have $D^p\in\Der(A)~\forall D\in \Der(A)$ and therefore $\Der(A)$ is a restricted Lie algebra with the bracket $[D,D']=D\circ D'-D'\circ D$ and the $p$-map $D^{[p]}=D^p.$
\end{Example}

\sssbegin{Remark}(See also \cite[Remark 1.1]{KL})\label{rmk:centerless}
    Suppose that $L$ is a centerless Lie algebra equipped with a map $(-)^{[p]}:L\rightarrow L$ satisfying Eq. \eqref{eq:def-2}. Then, $(-)^{[p]}$ automatically satisfies Eqs. \eqref{eq:def-1} and \eqref{eq:si}. Indeed, let $x,y\in L$. We have
    \begin{align*}
    \ad_{(x+y)^{[p]}}&=\ad^p_{x+y}=\ad^p_{x}+\ad^p_{y}+\sum_{i=1}^{p-1}\overline{s_i}(\ad_x,\ad_y)\\
            &=\ad_{x^{[p]}}+\ad_{y^{[p]}}+\sum_{i=1}^{p-1}\ad_{s_i(x,y)},
    \end{align*} where $\overline{s_i}$ refers to the commutator bracket (note that the equality $\sum_{i=1}^{p-1}\overline{s_i}(\ad_x,\ad_y)=\sum_{i=1}^{p-1}\ad_{s_i(x,y)}$ follows from the Jacobi identity). It follows that $$(x+y)^{[p]}-x^{[p]}-y^{[p]}-\sum_{i=1}^{p-1}s_i(x,y)$$ lies in the center of $L$, which is trivial. 
\end{Remark}
\subsection{The Hochschild Lemma} The following general statement appeared in Hochschild's work on simple algebras with purely inseparable splitting fields of exponent 1, see \cite[Lemma 1]{Ho}.\\
\begin{center}
\begin{minipage}{0.85\textwidth}
Let $U$ be an associative algebra over the field of integers modulo $p$ and let $V\subset U$ be a commutative subalgebra. For all $u\in U$, denote by $D_u$ the derivation defined by $D_u(w)=uw-wu,~\forall w\in U$. Then, for all $u\in U$ such that $D_u(V)\subset V$, we have
$(vu)^p=v^pu^p+D_{vu}^{p-1}(v)u,~\forall v\in V.$\\
\end{minipage}
\end{center}
In the context of Lie-Rinehart algebras, this statement will be formulated as follows. Let $(A,L,\rho)$ be a Lie-Rinehart algebra. Applying the Lemma with $U$ the universal enveloping algebra of $(A,L,\rho)$ (see \cite{Ri} or Section \ref{sec:UEA} below), and $V=A$, we obtain
    \begin{equation}\label{eq:Hoch-old}
    (ax)^{p}=a^px^{p}+\rho(ax)^{p-1}(a)x,\quad\forall a\in A,~\forall x\in L.
\end{equation}
A modern version of Hochschild's Lemma appears then in Schauenburg's work as follows.

\sssbegin{Proposition}[\text{\cite[Lemma 2.1]{Sc}}]\label{prop:hochschild}
    Let $(A,L,\rho)$ be a Lie-Rinehart algebra and let $(\phi,M)$ be a Lie-Rinehart module. Then, we have
    \begin{equation}\label{eq:Hoch-new}
        \phi(ax)^p=a^p\phi(x)^p+\rho(ax)^{p-1}(a)\phi(x).
    \end{equation}
\end{Proposition} 

Equation \eqref{eq:Hoch-old} (and \eqref{eq:Hoch-new}) serves as the main motivation for the definition of restricted Lie-Rinehart algebras below.

\subsection{Restricted Lie-Rinehart algebras}\label{def:RLR}
Following \cite{Ru,Do}, a Lie-Rinehart algebra $(A,L,\rho)$ (see Eq. \eqref{eq:def-LR}) is called \textit{restricted} if $L$ is a restricted with a $p$-map $(-)^{[p]}:L\rightarrow L$, and if moreover, we have
\begin{equation}\label{eq:Hoch}
    (ax)^{[p]}=a^px^{[p]}+\rho(ax)^{p-1}(a)x,\quad\forall a\in A,~\forall x\in L.
\end{equation}

Equation \eqref{eq:Hoch} is referred to as \textit{Hochschild condition} and is derived from the Hochschild's Lemma \cite[Lemma 1]{Ho}.


\sssbegin{Example}[see \cite{Do}]
    Let $A$ be an associative commutative algebra and $\Der(A)$ its restricted Lie algebra of derivations. Then, $(A,\Der(A),\id)$ is a Lie-Rinehart algebra. Applying Hochschild's Lemma with $\phi=\rho=\id$ (see Eq. \eqref{eq:super-hochschild-ev/ev}), we have
    \begin{equation}
        (aD)^p=a^pD^p+(aD)^{p-1}(a)D,\quad \forall a\in A,~\forall D\in \Der(A).
    \end{equation}
    Therefore, $(A,\Der(A),\id)$ is a restricted Lie-Rinehart algebra.
\end{Example}

\sssbegin{Example}
    Let $(A,\{-,-\},)$ be a Poisson algebra with unit $1_A$ equipped with a $p$-map $(-)^{\{p\}}$ (that is, $(A,\{-,-\},(-)^{\{p\}})$ is a restricted Lie algebra). Following \cite{BYZ}, $(A,\{-,-\},(-)^{\{p\}})$ is called \textit{restricted Poisson algebra} if
    \begin{align}
        (x^2)^{\{p\}}&=2x^px^{\{p\}},\quad\forall x,y\in A.
    \end{align}
 Consider the module of Kähler differentials $\Omega^1(A)$ of a restricted Poisson algebra $(A,\{\cdot,\cdot\},(\cdot)^{\{p\}})$, which is the $A$-module generated by symbols $du,~u\in A,$ with relations
    \begin{equation}
        d(u+v)=du+dv,~d(uv)=udv+vdu,~d1_A=0.
    \end{equation}
    Then, it has been shown in \cite[Theorem 3.8]{Hu} that $(A,\Omega^1(A),\rho)$ is a Lie-Rinehart algebra, with
    \begin{align}
     [xdu,ydv]_{\Omega^1(A)}&:=x\{u,y\}d v+y\{x,v\}d u+xyd\{u,v\},\quad &\forall x,y,u,v\in A;\\
     \rho(xdu)&:=x\{u,-\} &\forall x,u\in A.
    \end{align}
    Suppose that $\Omega^1(A)$ is free as an $A$-module. Then, $(A,\Omega^1(A),\rho)$ becomes a restricted Lie-Rinehart algebra, with the $p$-map (see \cite[Theorem 8.2]{BYZ}) 
    \begin{equation}
        (xdu)^{[p]}:=x^pd(u^{\{p\}})+\rho(xdu)^{p-1}(x)du,\quad \forall x,u\in A.
    \end{equation}
    The analogue construction in the case where $p=2$ was carried out  in \cite{BEL}.
\end{Example}

\section{A superization of Hochschild's Lemma}\label{sec:super-hochschild}
In this Section, we give a superized version of Hochschild's Lemma (see \cite[Lemma 1]{Ho}, \cite[Lemma 2.1]{Sc}, also Proposition \ref{prop:hochschild}). Our approach follows Schauenburg's proof (see Section 5 of \cite{Sc}). 

\ssbegin{Theorem}[Superization of the Hochschild's Lemma] \label{thm:super-hochschild} Let $(A,L,\rho)$ be a Lie-Rinehart superalgebra over a field $\K$ of characteristic $p\geq 3$ and let $(\phi,M)$ be a Lie-Rinehart module. Then, we have 
\begin{align}\label{eq:super-hochschild-ev/ev}
    \phi(ax)^p&=a^p\phi(x)^p+\rho(ax)^{p-1}(a)\phi(x),&\forall a\in A_\ev,~\forall x\in L_\ev;\\
    \phi(ax)^{2p}&=a^{2p}\phi(x)^{2p}+\rho(ax)^{2p-1}(a)\phi(x)&\\\nonumber
    &~~+\sum_{i=0}^{p-1}\lambda_{i}\rho(ax)^i(a)\rho(ax)^{2p-2-i}(a)\phi(x)^{2},&\forall a\in A_\ev,~\forall x\in L_\od;\\
    \phi(ax)^{2p}&=0, &\quad\forall a\in A_\od,~\forall x\in L_\ev;\\
     \phi(ax)^p&=a\bigl(\rho(x)(a)\bigl)^{p-1}\phi(x),\quad&\forall a\in A_\od,~\forall x\in L_\od,
\end{align}
where the coefficients $\lambda_{i}$ are given by 
\begin{equation} \label{eq:lambda}
\lambda_{i}=
\begin{cases}
    2(-1)^{\frac{i}{2}}~&\text{ if $i$ is even, } 0\leq i<p-1; \\[1mm]
    2(-1)^{\frac{i-1}{2}}~&\text{ if $i$ is odd, } 1\leq i<p-1;\\[1mm]
    (-1)^{\frac{p-1}{2}} &\text{ if } i=p-1.
    \end{cases}
\end{equation}

\end{Theorem}

\sssbegin{proof}
    The purely even case follows from \cite[Lemma 2.1]{Sc}. The remaining cases are covered by Lemmas \ref{lem:even/odd}, \ref{lem:odd/even} and \ref{lem:odd/odd} below.
\end{proof}

The rest of this Section is devoted to the proof of \ref{thm:super-hochschild}. The study is divided into four different cases according to the parity of the elements involved.

\subsection{The general setting} \label{sec:general-setting} Let $(A,L,\rho)$ be a Lie-Rinehart superalgebra. Consider the $\Z_2$-graded ring $V=\Z[x_0,x_1,x_2,\cdots]$, such that $|x_{i+1}|=|x_i|+|\delta|$, where $\delta$ is the derivation of $V$ defined by $\delta(x_i)=x_{i+1}$.
Consider the superbialgebra $H=\Z[\delta]$ given by the coproduct $\Delta(1)=1\otimes 1,$ and $\Delta(\delta)=1\otimes\delta+\delta\otimes 1.$ Recall that according to the Koszul rule, the product on $H$ is given by \begin{equation}
(a\otimes b)(c\otimes d)=(-1)^{|b||c|} ac\otimes bd,\quad\forall~a,b,c,d\in H.
\end{equation} 
Recall also that the coproduct satisfies
\begin{equation}
\Delta(H_\ev)\subset H_\ev\otimes H_\ev+H_\od\otimes H_\od,\quad \Delta(H_\od)\subset H_\ev\otimes H_\od+H_\od\otimes H_\ev.
\end{equation} 

We can form the smash product $V\#H$ as follows. Using Sweedler's notation, the product in $V\#H$ is given by
    \begin{equation}(v\#g)(w\#h)=\sum_{(g)}(-1)^{|w||g_{(2)}|}v(g_{(1)}w)\#(g_{(2)}h),~\forall v,w\in V,\quad\forall g,h\in H.\end{equation} In particular, if $g$ is primitive (that is, $\Delta(g)=1\otimes g+g\otimes 1$), the above product becomes
    \begin{equation}(v\#g)(w\#h)=(-1)^{|w||g|}vw\#gh+v(gw)\#h,\quad\forall v,w\in V,~\forall g,h\in H.\end{equation}

Now, let $a\in A$ such that $|a|=|x_0|$ and let $\alpha\in\Der(A)$ such that $|\alpha|=|\delta|$. Consider the (unique) ring homomorphism $f:V\rightarrow A$ given by $f(x_0)=a$ satisfying $f\circ\delta=\alpha\circ f$. From this equality, we obtain $f(x_i)=\alpha^i(a),~\forall i\geq 0$.

Let $M$ be an $A$-module and let $\beta\in\End_{\Z}(M)$ such that $$\beta(am)=\alpha(a)m+(-1)^{|a||\beta|}a\beta(m),\quad\forall a\in A,~\forall m\in M.$$
Note that such maps $\beta$ do exist, for example, if $(\phi,M)$ is a representation of the Lie-Rinehart superalgebra $(A,L,\rho)$ and $x\in L$ homogeneous, take $\beta=\phi(x)$ and $\alpha=\rho(x)$.

\sssbegin{Lemma}\label{lem:superaction}
The $A$-module $M$ becomes a $V\#H$-module with \begin{equation}
            (v\#1)m=f(v)m,\quad \text{and} \quad  (1\#\delta)m=\beta(m),\quad \forall v\in V,~\forall m\in M.
        \end{equation}
\end{Lemma}

\begin{proof}
    Let $v,w\in V$ and
    $m\in M$. We have $$(v\#1)(w\#1)m=(v\#1)f(w)m=f(v)f(w)m=f(vw)m=((vw)\#1)m.$$  The module structure is then extended by
        $$(1\#\delta^k)m=\beta^k(m)~\forall k\geq 1,\quad\text{and}\quad (v\#\delta)m=(v\#1)(1\#\delta)m.$$
\end{proof}

\sssbegin{Lemma}\label{lem:recursive-Gamma}
            We have $(x_0\#\delta)^k=\sum_{j=1}^k\Gamma_{k,j}\#\delta^j,$ where $\Gamma_{k,j}\in\Z[x_0,x_1,\cdots, x_k]$ are recursively defined by $\Gamma_{1,1}=x_0$ and
           $$ \Gamma_{k+1,j}=\begin{cases}
            x_0\delta(\Gamma_{k,1}),&j=1\\[1mm]
            x_0\delta(\Gamma_{k,j})+(-1)^{|\delta||\Gamma_{k,j-1}|}x_0\Gamma_{k,j-1},&j=2,\cdots,k\\[1mm]
            x_0\Gamma_{k,k},&j=k+1.
            \end{cases}$$ In particular, we have $\Gamma_{k,k}=x_0^k$ and $\Gamma_{k,1}=(x_0\delta)^{k-1}(x_0).$
    \end{Lemma}
    \begin{proof}
We will use induction. First, notice that the elements $\Gamma_{k,j}$ are homogeneous of parity $|\Gamma_{k,j}|=k-j~\mod{2}.$ The induction  step is given by
       \begin{align*}
(x_0\#\delta)\sum_{j=1}^k\Gamma_{k,j}\#\delta^j&=\sum_{j=1}^k(-1)^{|\delta||\Gamma_{k,j}|}x_0\Gamma_{k,j}\#\delta^{j+1}+x_0\delta(\Gamma_{k,j})\#\delta^j\\
&=x_0\delta(\Gamma_{k,1})\#\delta+\sum_{j=2}^k\Bigl((-1)^{|\delta||\Gamma_{k,j-1}|}x_0\Gamma_{k,j-1}\#\delta^j+x_0\delta(\Gamma_{k,j})\#\delta^j\Bigl)\\
          &+(-1)^{|\delta||\Gamma_{k,k}|}x_0\Gamma_{k,k}\#\delta^{k+1}.
       \end{align*}
The result follows.     
    \end{proof}

\sssbegin{Corollary}
For all $k\geq 2$, we have
\begin{equation}
    \Gamma_{k,2}=(-1)^{|\delta||\Gamma_{k-1,1}|}x_0(x_0\delta)^{k-2}(x_0)+x_0\delta(\Gamma_{k-1,2}).
\end{equation}
\end{Corollary}


\subsection{The case where $x$ is odd and  $a$ is even} Let $(A,L,\rho)$ be a Lie-Rinehart superalgebra. In this section, we investigate the case where $a\in A$ is even and $x\in L$ is odd.

\sssbegin{Proposition}\label{prop:p-divides-even/odd}
    Let $p\geq 3$ be a prime number. Suppose that $x_0$ is even and that $\delta$ is odd. Then, $\Gamma_{2p,j}\equiv 0\mod p,$ for $3\leq j\leq 2p-1,~j\neq p$.
\end{Proposition}
\begin{proof} We use the construction of Section \ref{sec:general-setting} by choosing $M=A=V$ and $\alpha=\beta=\delta$. In that case, $f=\id$.
        First, notice that since $\delta$ is an odd derivation, then $\delta^2$ is an even derivation. Moreover, one can show that $(x_0\#\delta)$ acts on $M$ as an odd derivation and consequently $(x_0\#\delta)^{2p}$ is an even derivation. Let $m>2p$ such that $x_m\in V$ is even. For $x_r,x_s\in V$ and $j\geq 1$, we have
        $$x_r\bigl(\Gamma_{2p,j}\#\delta^j\bigl)(x_s)=x_r\bigl(\Gamma_{2p,j}\#1\bigl)\delta^j(x_s)=x_r\Gamma_{2p,j}\delta^j(x_s)=(-1)^{|x_r||\Gamma_{2p,j}|}\Gamma_{2p,j}x_r\delta^j(x_s).$$

        It follows that
        \begin{align*}
            0&=(x_0\#\delta)^{2p}(x_0x_m)-x_0(x_0\#\delta)^{2p}(x_m)-x_m(x_0\#\delta)^{2p}(x_0)\\
            &=\sum_{j=1}^{2p}\Bigl(\bigl(\Gamma_{2p,j}\#\delta^j\bigl)(x_0x_m)-x_0\bigl(\Gamma_{2p,j}\#\delta^j\bigl)(x_m)-x_m\bigl(\Gamma_{2p,j}\#\delta^j\bigl)(x_0)\Bigl)\\
            &=\sum_{j=1}^{2p}\Gamma_{2p,j}\bigl(\delta^j(x_0x_m)-x_0\delta^j(x_m)-x_m\delta^j(x_0)\bigl)
            \\
            &=\sum_{j=3}^{2p-1}\Gamma_{2p,j}\bigl(\delta^j(x_0x_m)-x_0\delta^j(x_m)-x_m\delta^j(x_0)\bigl).\\
            &=\sum_{j=3}^{2p-1}\Gamma_{2p,j}\sum_{i=1}^{j-1}\binom{j}{i}\delta^i(x_0)\delta^{j-i}(x_m)
        \end{align*} 
        
        (Note that since $\delta, \delta^{2}$ and $\delta^{2p}$ are  derivations, the terms for $j=1,2$ and $j=2p$ vanish.)
         Let $\ell\in\{3,\cdots,2p-1\}$ such that $\Gamma_{2p,j}=0$ for all $3\leq j\leq \ell-1,~j\neq p$. Then, we have
        \begin{align*}0&=\Gamma_{2p,\ell}\sum_{i=1}^{\ell-1}\binom{\ell}{i}\delta^i(x_0)\delta^{\ell-i}(x_m)+\sum_{j=\ell+1}^{2p-1}\Gamma_{2p,j}\sum_{i=1}^{j-1}\binom{j}{i}\delta^i(x_0)\delta^{j-i}(x_m)\\
        &=\Gamma_{2p,\ell}\binom{\ell}{1}x_1x_{m+\ell-1}+P,
        \end{align*} with $P$ some polynomial that does not contain the monomial $x_1x_{m+\ell-1}$. Therefore, we have $\Gamma_{2p,j}\equiv 0\mod p,$ for $3\leq j\leq 2p-1,~j\neq p$.       
\end{proof}
In the case where $j=p$, the above proof no longer works, so a different approach is required. This case is addressed in the following Proposition; its proof is deferred to Appendix~\ref{sec:gamma2pp}.

\sssbegin{Proposition} \label{prop:p-divides-Gamma2pp}
 Let $p\geq 3$ be a prime number. Suppose that $x_0$ is even and that $\delta$ is odd. Then, $\Gamma_{2p,p}\equiv 0\mod p.$
\end{Proposition}
\begin{proof}
    See Appendix \ref{sec:gamma2pp}.
\end{proof}

Using Propositions \ref{prop:p-divides-even/odd} and \ref{prop:p-divides-Gamma2pp}, we have
\begin{equation}
    (x_0\#\delta)^{2p}=\Gamma_{2p,1}\#\delta+\Gamma_{2p,2}\#\delta^2+\Gamma_{2p,2p}\#\delta^{2p}.
\end{equation} Let $(\phi,M)$ be a representation of $(A,L,\rho)$, let $m\in M$ and let $x\in L_\od$. Using Lemma \ref{lem:superaction} with $\beta=\phi(x)$ and $\alpha=\rho(x)\in\Der(A)$, we have
\begin{equation}\label{eq:2p-even/odd}
    (x_0\#\delta)^{2p}(m)=f\bigl(\Gamma_{2p,1}\bigl)\phi(x)(m)+f\bigl(\Gamma_{2p,2}\bigl)\phi(x)^2(m)+f\bigl(\Gamma_{2p,2p}\bigl)\phi(x)^{2p}(m).
\end{equation} Therefore, it remains to compute the terms $f\bigl(\Gamma_{2p,j}\bigl)$ for $j\in\{1,2,2p\}$.

\sssbegin{Lemma}\label{lem:even/od-A}
        Let $(A,L,\rho)$ be a Lie-Rinehart superalgebra. Let $a\in A_\ev,~x\in L_\od$ and $k\geq 1$. Then, we have
        \begin{equation}
            f\bigl(\Gamma_{k,1}\bigl)=\rho(ax)^{k-1}(a),\quad\text{and}\quad f\bigl(\Gamma_{k,k}\bigl)=a^k. 
        \end{equation}
\end{Lemma}

\begin{proof} 
Let $v\in V$. We have
$$f\bigl((x_0\delta)(v)\bigl)=f(x_0)f\circ\delta(v)=a\rho(x)\circ f(v).$$ Applying the relation to $v=(x_0\delta)^{k-2}(x_0)$, we obtain by direct induction
\begin{align*}
    f\bigl(\Gamma_{k,1}\bigl)=f\bigl((x_0\delta)^{k-1}(x_0)\bigl)&=a\rho(x)f\bigl((x_0\delta)^{k-2}(x_0)\bigl)\\
    &=(a\rho(x))^{k-1}f(x_0)=\rho(ax)^{k-1}(a).
\end{align*} Moreover, we have $f\bigl(\Gamma_{k,k}\bigl)=f(x_0^k)=a^k.$
\end{proof}

\sssbegin{Lemma}\label{lem:even/od-B} Let $(A,L,\rho)$ be a Lie-Rinehart superalgebra and let $a\in A_\ev,$ and $x\in L_\od$. Then, we have $f(\Gamma_{3,2})=a\rho(ax)(a)$ and for $k\geq 4$, we have
\begin{equation}\label{eq:x-odd-a-even}
f(\Gamma_{k,2})=\sum_{i=0}^{k-2}\mu_{k,i}\rho(ax)^i(a)\rho(ax)^{k-2-i}(a),
\end{equation}
where $\mu_{3,0}=\mu_{3,1}=\frac{1}{2}$ and for $k\geq 4$,
\begin{equation}\mu_{k,0}=\mu_{k-1,0}+(-1)^{k};\quad \mu_{k,i}=\mu_{k-1,i-1}+(-1)^i\mu_{k-1,i},~ 1\leq i\leq k-3;\quad \mu_{k,k-2}=\frac{1}{2}\end{equation}
\end{Lemma}
A direct calculation, for instance, gives after simplification:
$$
\begin{array}{lcl}
f(\Gamma_{5,2})&=&a\rho(ax)^3(a)+2\rho(ax)(a)\rho(ax)^2(a),\\[1mm]
f(\Gamma_{6,2}) & = & 2a\rho(ax)^4(a)-\rho(ax)(a)\rho(ax)^3(a)+2\rho(ax)^2(a)\rho(ax)^2(a),\\[1mm]
f(\Gamma_{10,2}) & = & 2 a\rho(ax)^8(a)-3 \rho(ax)(a)\rho(ax)^7(a) + 8 \rho(ax)^2(a)\rho(ax)^6(a)-2  \rho(ax)^3(a)\rho(ax)^5(a)\\[2mm]
&&+ 6 \rho(ax)^4(a)\rho(ax)^4(a).
\end{array}
$$
\begin{proof}
    Let $x\in L_\od, a\in A_\ev.$ Recall that $f\circ\delta=\alpha\circ f=\rho(x)\circ f$. Let $k\geq 4$ and suppose that we have
    $$f(\Gamma_{k,2})=\sum_{i=0}^{k-2}\mu_{k,i}\rho(ax)^i(a)\rho(ax)^{k-2-i}(a).$$ Then, we have
    \begin{align*}
        af\circ\delta(\Gamma_{k,2})&=\rho(ax)\Bigl(\sum_{i=0}^{k-2}\mu_{k,i}\rho(ax)^i(a)\rho(ax)^{k-2-i}(a)\Bigl)\\
        &=\sum_{i=0}^{k-2}\mu_{k,i}\rho(ax)^{i+1}(a)\rho(ax)^{k-2-i}(a)+\sum_{i=0}^{k-2}(-1)^i\mu_{k,i}\rho(ax)^i(a)\rho(ax)^{k-1-i}(a)\\
        &=\sum_{j=1}^{k-1}\mu_{k,j-1}\rho(ax)^{j}(a)\rho(ax)^{k-1-j}(a)+\sum_{i=0}^{k-2}(-1)^i\mu_{k,i}\rho(ax)^i(a)\rho(ax)^{k-1-i}(a)\\
        &=\mu_{k,0}a\rho(ax)^{k-1}(a)+\sum_{i=1}^{k-2}\bigl(\mu_{k,i-1}+(-1)^{i}\mu_{k,i}\bigl)\rho(ax)^i(a)\rho(ax)^{k-1-i}(a)\\
        &~\quad+\mu_{k,k-2}\rho(ax)^{k-1}(a)a.
    \end{align*} Therefore, we have
    \begin{align*}
        f(\Gamma_{k+1,2})=&~f\bigl((-1)^{|\delta|(k+1)}x_0(x_0\delta)^{k-1}+x_0\delta(\Gamma_{k,2})\bigl)\\
        =&~(-1)^{|\delta|(k+1)}a\rho(ax)^{k-1}(a)+af\circ\delta(\Gamma_{k,2})\\
        =&~\bigl(\mu_{k,0}+(-1)^{k+1}\bigl)a\rho(ax)^{k-1}(a)+\sum_{i=1}^{k-2}\bigl(\mu_{k,i-1}+(-1)^i\mu_{k,i}\bigl)\rho(ax)^i(a)\rho(ax)^{k-1-i}(a)\\
        &+\mu_{k,k-2}\rho(ax)^{k-1}(a)a,
    \end{align*} thus the conclusion.
\end{proof}
\sssbegin{Remark}\label{rmk:reduced-coefs}
    Since $A$ is a supercommutative superalgebra, Eq. \eqref{eq:x-odd-a-even} can be simplified. Indeed, suppose that $k=2p$ for $p\geq 3$ prime (this is the case that interests us). Then we have
    \begin{align*}
        f(\Gamma_{2p,2})&=\sum_{i=0}^{2p-2}\mu_{2p,i}\rho(ax)^i(a)\rho(ax)^{2p-2-i}(a)\\
        &=\sum_{i=0}^{p-2}\bigl(\mu_{2p,i}+(-1)^i\mu_{2p,2p-2-i}\bigl)\rho(ax)^i(a)\rho(ax)^{2p-2-i}(a)+\mu_{2p,p-1}\rho(ax)^{p-1}(a)^2\\
        &=\sum_{i=0}^{p-1}\lambda_{i}\rho(ax)^i(a)\rho(ax)^{2p-2-i}(a),
    \end{align*} with $\lambda_{i}:=\mu_{2p,i}+(-1)^i\mu_{2p,2p-2-i}, \text{ for } i\neq p-1$ and $\lambda_{p-1}:=\mu_{2p,p-1}$ the \textit{reduced coefficients}.
\end{Remark}

\sssbegin{Lemma} The reduced coefficients modulo $p$ satisfy
\begin{equation} \label{eq:lambda-bis}
\lambda_{i}=
\begin{cases}
    2(-1)^{\frac{i}{2}}~&\text{ if $i$ is even, } 0\leq i<p-1; \\[1mm]
    2(-1)^{\frac{i-1}{2}}~&\text{ if $i$ is odd, } 1\leq i<p-1;\\[1mm]
    (-1)^{\frac{p-1}{2}} &\text{ if } i=p-1.
    \end{cases}
\end{equation}
\end{Lemma}
\begin{proof}
    The coefficients $(\mu_{k,i})$ defined in Lemma \ref{lem:even/od-B} satisfy the relations
    \begin{align*}
        \mu_{2r+2,2j}&=\frac{1}{2}\binom{r}{j}+\binom{r-1}{j},\quad 0\leq j\leq r;& 
        \mu_{2r+2,2j+1}&=-\binom{r-1}{j+1},\quad 0\leq j\leq r-1;\\
        \mu_{2r+3,2j}&=\frac{1}{2}\binom{r}{j},\quad 0\leq j\leq r;& 
        \mu_{2r+3,2j+1}
    &=\frac{1}{2}\binom{r}{j}+\binom{r}{j+1},\quad 0\leq j \leq r.
    \end{align*}
    Recall that the following identities hold modulo $p$:
    \begin{equation*}
        \binom{p-1}{j}\equiv(-1)^j;\quad \text{and }  \binom{p-2}{j}\equiv(j+1)(-1)^j.
    \end{equation*}
   It follows that modulo $p$ we have 
    \begin{align}\label{eq:modp-lemma}
        \mu_{2p,2j}\equiv\bigl(\frac{1}{2}+j+1\bigl)(-1)^j;\quad \text{and } \mu_{2p,2j+1}\equiv (j+2)(-1)^j.
    \end{align} Let $k=2p-2-i$ and suppose that $i\neq p-1$.
    
    \underline{The case where $i$ is even}. Let $i=2j$ which implies that  $k=2(p-1-j).$ It follows from Eq. \eqref{eq:modp-lemma} that
    \begin{equation*}    
        \mu_{2p,k}=(-1)^{p-1-j}\bigl(\frac{1}{2}+p-1-j+1\bigl)\equiv(-1)^j\bigl(\frac{1}{2}-j\bigl)\mod p.
    \end{equation*} Thus, we have
    \begin{equation*}
       \displaystyle  \lambda_i=\mu_{2p,i}+(-1)^i\mu_{2p,k}\equiv 2(-1)^j=2(-1)^{\frac{i}{2}}\mod p.
    \end{equation*}
    
    \underline{The case where $i$ is odd}. Let $i=2j+1$. In that case, we have $k=2(p-2-j)+1$. Let $t=p-2-j$. We have
    \begin{equation*}
        \mu_{2p,k}=\mu_{2p,2t+1}=(t+2)(-1)^t=(p-j)(-1)^{p-2-j}\equiv j(-1)^j\mod p.
    \end{equation*} It follows that
      \begin{equation*}
       \lambda_i=(j+2)(-1)^j-j(-1)^j\equiv 2(-1)^j=2(-1)^{\frac{i-1}{2}}\mod p.\qed
    \end{equation*}
    \noqed
\end{proof}

\sssbegin{Example} The reduced coefficients $\lambda_i$ modulo $p$ for small values of $p$ are given in the following table.
\begin{center}
    \begin{tabular}{|c|c|c|c|}
    \hline
          &$p=3$&$p=5$&$p=7$\\\hline    
          $\lambda_{0}$&2&2&2\\ 
          $\lambda_{1}$&2&2&2\\
          $\lambda_{2}$&2&3&5\\
          $\lambda_{3}$&-&3&5\\
          $\lambda_{4}$&-&1&2\\
          $\lambda_{5}$&-&-&2\\
          $\lambda_{6}$&-&-&6\\\hline            
    \end{tabular}
   \\~\\ Reduced coefficients $\lambda_{i}$ modulo $p$ for small $p$'s.
\end{center}

\end{Example}

\sssbegin{Lemma}\label{lem:even/odd}
    Let $(A,L,\rho)$ be a Lie-Rinehart superalgebra and let $(\phi,M)$ be a Lie-Rinehart module. Then, for all $x\in L_\od$ and all $a\in A_\ev$, we have
    \begin{equation}
        \phi(ax)^{2p}=a^{2p}\phi(x)^{2p}+\sum_{i=0}^{p-1}\lambda_{i}\rho(ax)^i(a)\rho(ax)^{2p-2-i}(a)\phi(x)^{2}+\rho(ax)^{2p-1}(a)\phi(x),
    \end{equation}
    where the coefficients $\lambda_{i}$ are defined in Eq. \eqref{eq:lambda}.
\end{Lemma}

\begin{proof}
Let $x\in L_\od$ and all $a\in A_\ev$. Rewriting Eq. \ref{eq:2p-even/odd} using Lemmas \ref{lem:even/od-A} and \ref{lem:even/od-B}, we have for all $m\in M$
\begin{align}
    (x_0\#\delta)^{2p}(m)=\rho&(ax)^{2p-1}(a)\phi(x)(m)+\sum_{i=0}^{p-1}\lambda_{i}\rho(ax)^i(a)\rho(ax)^{2p-2-i}(a)\phi(x)^2(m)\\\nonumber&+a^{2p}\phi(x)^{2p}(m).
\end{align}
Moreover, a direct computation yields $(x_0\#\delta)^{2p}(m)=\phi(ax)^{2p}(m),\quad\forall m\in M.$
\end{proof}
\subsection{The case where $x$ is even and  $a$ is odd}
Let $(A,L,\rho)$ be a Lie-Rinehart superalgebra. In this section, we investigate the case where $a\in A$ is odd and $x\in L$ is even. In that case, $V=\Z[x_0,x_1,x_2,\cdots]$ with all $x_i$'s being odd.
\sssbegin{Lemma}\label{lem:odd/even}
    Let $(A,L,\rho)$ be a Lie-Rinehart superalgebra and let $(\phi,M)$ be a Lie-Rinehart module. Then, we have
    \begin{equation}
    \phi(ax)^{2p}=0, \quad\forall x\in L_\ev,~\forall a\in A_\od.
    \end{equation}
\end{Lemma}

\begin{proof}
Let $\delta$ be the even derivation sending $x_i$ to $x_{i+1}$. The Proposition follows from the fact that $\Gamma_{k,k}=0$ for all $k\geq 1$ and $\Gamma_{k,j}=0$ for all $k\geq 3$ and all $j=1,\cdots, k$. Indeed, we have $\Gamma_{1,1}=x_0$, $\Gamma_{2,1}=x_0x_1$ and $\Gamma_{k,k}=x_0^k=0,~\forall k\geq 2$. It follows that $\Gamma_{3,1}=x_0x_1^2+x_0^2x_1=0,$ since $x_0$ and $x_1$ are odd. Consequently, $\Gamma_{k,1}=0,~\forall k\geq 3$. Moreover, we have
 $$\Gamma_{3,2}=x_0\delta(\Gamma_{2,2})+x_0\Gamma_{2,1}=0.$$
By induction, it follows that $\Gamma_{k,j}=0$ for all $k\geq 3$ and all $j=1,\cdots, k$. In particular, all the terms of the form $\Gamma_{2p,j}$ (for $p\geq 3$ prime) vanish.
\end{proof}

\subsection{The case where $x$ and $a$ are odd}
Let $(A,L,\rho)$ be a Lie-Rinehart superalgebra. In this section, we investigate the case where $a\in A$ and $x\in L$ are odd. In that case, $V=\Z[x_0,x_1,x_2,\cdots]$ with $|x_0|=1$ and $|x_i|=i+1\mod2.$

\sssbegin{Lemma}\label{lem:odd/odd}
      Let $(A,L,\rho)$ be a Lie-Rinehart superalgebra and let $(\phi,M)$ be a Lie-Rinehart module. Then, we have
      \begin{equation}
            \phi(ax)^p=a\bigl(\rho(x)(a)\bigl)^{p-1}\phi(x),\quad\forall a\in A_\od,~\forall x\in L_\od.
      \end{equation}
\end{Lemma}

\begin{proof}
    Let $\delta$ be the odd derivation sending $x_i$ to $x_{i+1}$. Then, we have $\Gamma_{k,1}=x_0x_1^{k-1},~\forall k\geq 1,$ and $|\Gamma_{k,1}|=\od.$ 
    We claim that $\Gamma_{k,j}=0$ for $2\leq j\leq k$. Indeed, we have $\Gamma_{k,k}=x_0^k=0~\forall k\geq2$ and therefore $\Gamma_{2,2}=0.$ Next, we have $x_0\Gamma_{2,1}=x_0^2x_1=0.$ For the recursive step, suppose that $\Gamma_{k,j}=x_0\Gamma_{k,j-1}=0,~\text{for}~ 2\leq j\leq k$. Using Lemma \ref{lem:recursive-Gamma}, we obtain
    $$ \Gamma_{k+1,j}=x_0\delta(\Gamma_{k,j})-x_0\Gamma_{k,j-1}=0.$$ The claim is proved.
    Using Lemma \ref{lem:recursive-Gamma}, we have  $$(x_0\#\delta)^p=\sum_{j=1}^p\Gamma_{p,j}\#\delta^j=\Gamma_{p,1}\#\delta=x_0x_1^{p-1}\#\delta. $$ It follows that
    $$\phi(ax)^p(m)=(x_0\#\delta)^pm=(x_0x_1^{p-1}\#\delta)m=a\bigl(\rho(x)(a)\bigl)^{p-1}\phi(x)m,\quad\forall m\in M.$$
\end{proof}

\section{Restricted Lie-Rinehart superalgebras}\label{sec:sRLR}
We are now in position to define restricted Lie-Rinehart superalgebras. We start by recalling some basic facts about restricted Lie superalgebras.

\subsection{Restricted Lie superalgebras}

Let $L$ be a~finite-dimensional  Lie superalgebra defined over a field of characteristic $p>2$. For $p=3$, the Jacobi identity does not imply $[x,[x,x]]=0$, for all $x\in L_\od$, therefore we require it as part of the definition. 

Following \cite{Mi,Pe}, we say that $L$ has a~\textit{$p|2p$-structure} if $L_\ev$ is  a restricted Lie algebra with a $p$-map $(-)^{[p]}$ and 
\begin{equation}
\label{RRRS}
\ad_{x^{[p]}}(y)=(\ad_x)^p(y),\quad\forall x \in L_\ev, ~\forall y\in L.
\end{equation}
We set 
\[
(-)^{[2p]}:L_\od \rightarrow L_\ev, \quad x\mapsto (x^2)^{[p]},\quad  \text{ where } x^2:=\frac{1}{2}[x,x].
\]
The pair $\bigl(L, (-)^{[p|2p]}\bigl)$ is called a~\textit{restricted} Lie superalgebra. When the base field has characteristic $p=2$, there are several notions of restrictedness (see \cite{BLLS}); however, these will not be considered in the present work. For a description of $p|2p$-structures on simplie Lie superalgebras, see \cite{BKLLS}.

The following theorem is a~straightforward superization of Jacobson's theorem \ref{Jac}.
\sssbegin{Theorem}\label{SJac}
Let $(e_j)_{j\in J}$ be a~basis of $L_\ev$, and let the elements $f_j\in L_\ev$ be such that  ${(\ad_{e_j})^p=\ad_{f_j}}$. Then, there exists exactly one $p|2p$-mapping $(-)^{[p|2p]}:L\rightarrow L$ such that 
\[
e_j^{[p]}=f_j \quad \text{ for all $j\in J$}.
\]
\end{Theorem}

Let $\bigl(L, (-)^{[p|2p]_L}\bigl)$ and $\bigl(H, (-)^{[p|2p]_H}\bigl)$ be two restricted Lie superalgebras. A morphism of Lie superalgebras $\varphi: L\rightarrow H$ is called \textit{restricted} if in addition to $\varphi\bigl([x,y]\bigl)=\bigl[\varphi(x),\varphi(y)\bigl]~\forall x,y\in L$, we also have 
$$\varphi\bigl(x^{[p]_L}\bigl)=\varphi(x)^{[p]_H},~\forall x\in L_\ev.  $$
A homogeneous ideal $I=I_\ev\oplus I_\od$ of $L$ is called a~\textit{$p$-ideal} if it is closed under the $p|2p$-map; namely 
\[
\text{$x^{[p]}\in I_\ev$ for all $x\in I_\ev$}. 
\]
As $x^{[2p]}=(x^2)^{[p]}$, for all $x\in L_\od$, this would imply that 
\[
x^{[2p]}\in I_\ev \text{ for all $x\in I_\od$}.
\]
The \textit{center} of a restricted Lie superalgebra $L$ is the $p$-ideal 
$$\fz(L)=\bigl\{ x\in L,~ [x,y]=0~\forall y\in L \bigl\}.   $$

Let $L$ be a restricted Lie superalgebra. An $L$-module $M$ is called \textit{restricted} if
\[
\underbrace{x\cdots x}_{p\text{~~times}} \cdot m =x^{[p]}\cdot m \quad \text{for all $x\in L_\ev$ and $m\in M$.}
\]
\subsection{Restricted Lie-Rinehart superalgebras}\label{def:sRLR} Based on Theorem \ref{thm:super-hochschild}, we introduce the notion of a restricted Lie-Rinehart superalgebra as follows. A Lie-Rinehart superalgebra $(A,L,\rho)$ is called \textit{restricted} if the Lie algebra $(L,[-,-],(-)^{[p]})$ is restricted and if in addition, we have 
\begin{align}\label{eq:superhochschild-def-even/even}
    (ax)^{[p]}&=a^px^{[p]}+\rho(ax)^{p-1}(a)x,&\forall a\in A_\ev,~\forall x\in L_\ev;\\\label{eq:superhochschild-def-even/odd}
    (ax)^{[2p]}&=a^{2p}x^{[2p]}+\rho(ax)^{2p-1}(a)x&\\\nonumber
    &~~+\sum_{i=0}^{p-1}\lambda_{i}\rho(ax)^i(a)\rho(ax)^{2p-2-i}(a)x^{2},&\forall a\in A_\ev,~\forall x\in L_\od;
    \\
    \label{eq:superhochschild-def-odd/even}
    (ax)^{[2p]}&=0, &\quad\forall a\in A_\od,~\forall x\in L_\ev;\\\label{eq:superhochschild-def-odd/odd}
     (ax)^{[p]}&=a\bigl(\rho(x)(a)\bigl)^{p-1}x,\quad&\forall a\in A_\od,~\forall x\in L_\od;
    \end{align}
    where the coefficients $\lambda_{i}$ are defined in Eq. \eqref{eq:lambda}, and $x^2:=\frac{1}{2}[x,x]$.
In the case where $A$ and $L$ are purely even, we recover Definition \ref{def:RLR}.

\sssbegin{Example}\label{ex:der}
    Let $A$ be a supercommutative associative superalgebra and consider the Lie-Rinehart superalgebra $(A,\Der(A),\id)$. The Lie superalgebra $\Der(A)$ is actually restricted, and it follows from Theorem \ref{thm:super-hochschild} applied to the representation $(\id,A)$ that $(A,\Der(A),\id)$ is a restricted Lie-Rinehart superalgebra.

Let us detail the computation in the case where $p=3$. Let $a\in A_\ev$ and $D\in \Der(A)_\od$. A long but direct computation gives
    \begin{equation*}
        (aD)^6=a^6D^6+(aD)^5(a)D+2a^4D^3(a)D(a)D^2+2a^5D^4(a)D^2.
    \end{equation*} However, one can show that
    \begin{align*}
        2a(aD)^4(a)D^2&=a^4D^2(a)^2D^2+2a^4D(a)D^3(a)D^2+2a^5D^4(a)D^2;\\
        2(aD)(a)(aD)^3(a)D^2&=a^4D^3(a)D(a)D^2;\\
        2(aD)^2(a)^2D^2&=2a^4D^2(a)^2D^2.
    \end{align*} It follows that
    \begin{align*}
    2a^4D^3(a)D(a)D^2+2a^5D^4(a)D^2=2&a(aD)^4(a)D^2+2(aD)(a)(aD)^3(a)D^2\\&+2(aD)^2(a)^2D^2,
    \end{align*} and therefore we obtain
     \begin{equation*}
        (aD)^6=a^6D^6+(aD)^5(a)D+2a(aD)^4(a)D^2+2(aD)(a(aD)^3(a)D^2+2(aD)^2(a)^2D^2.
    \end{equation*}
   We therefore recover Eq. \eqref{eq:superhochschild-def-even/odd}.
Similarly, for $a\in A_{\bar 1}$ and $D\in \Der(A)_{\bar 1}$, we get 
\[
(aD)^3=(aD)(aD(aD))=aD(a)^2D,
\]
and we recover Eq. \ref{eq:superhochschild-def-odd/odd}. Furthermore, Eq. \eqref{eq:superhochschild-def-odd/even} is straightforward and for Eq. \eqref{eq:superhochschild-def-even/even} we refer to \cite{Do, Sc}.
    \end{Example}
    \sssbegin{Example}[The restricted Witt superalgebra] Let $n\geq 1$. Recall that the free commutative superalgebra on $n$ odd generators $\xi_1,\cdots, \xi_n$, denoted by $\Lambda(n)$ is isomorphic to the Grassmann algebra. Following \cite[Section 1.2]{SZ}, the restricted Witt superalgebra of rank $n$ denoted by $W(n)$ is defined as the Lie superalgebra of derivations of $\Lambda(n)$. The $p$-map on the even part $W(n)_\ev$ is given by the usual $p$-th power of derivations. It follows from Example \ref{ex:der} that $\bigl(\Lambda(n),W(n),\id\bigl)$ is a restricted Lie-Rinehart superalgebra.
    \end{Example}
    
Following are two examples adapted from \cite{EM1}.
     \sssbegin{Example}\label{ex:1121}
Consider the $(1|1)$-dimensional associative supercommutative superalgebra $A$ spanned by elements $e_1|e_2$ (even$|$odd), where $e_1$ is the unit and $e_2e_2=0.$ Consider the $2|1$-dimensional restricted Lie superalgebra $L$ spanned by elements $x_1,x_2|x_3$ equipped with the brackets and a $p$-map as follows: 
\[
[x_1,x_3]=x_3,~[x_2,x_3]=-x_3, \;\text{ and } x_1^{[p]}=x_1+\alpha(x_1+x_2),~x_2^{[p]}=-x_1+\alpha(x_1+x_2),
\]
where $\alpha\in \K$. Then, $L$ is an $A$-module with the action
$$e_1x=x~\forall x\in L,\quad e_2x_1=\beta x_3,~\beta\in\K,\quad e_2x_2=\gamma x_3,~\gamma\in\K.$$ The pair $(A,L,\rho)$ is a restricted Lie-Rinehart superalgebra, with the anchor $\rho$ given by
$$\rho(x_1)(e_2)=e_2,\quad \rho(x_2)(e_2)=-e_2.$$
\end{Example}

     \sssbegin{Example}
Consider the $(1|1)$-dimensional associative supercommutative superalgebra $A$ of Example \ref{ex:1121}, and consider the $2|2$-dimensional restricted Lie superalgebra $L$ spanned by elements $x_1,x_2|x_3,x_4$ equipped with the brackets and a $p$-map as follows: \[
[x_1,x_3]=x_3,~[x_1,x_4]=x_4~[x_2,x_4]=x_3,\; \text{ and } x_1^{[p]}=x_1,~x_2^{[p]}=0.
\]
Then, $L$ is an $A$-module with the action
$$e_1x=x~\forall x\in L,\quad e_2x_1=\alpha x_3,~\alpha\in\K,\quad e_2x_2=\beta x_3,~\beta\in\K.$$ The pair $(A,L,\rho)$ is a restricted Lie-Rinehart superalgebra, with the anchor $\rho$ given by
$$\rho(x_1)(e_2)=e_2.$$
    \end{Example}
\black

\subsection{Modules and representations}\label{def:rest-rep} Let $(A,L,\rho)$ be a restricted Lie-Rinehart superalgebra. A representation of $(A,L,\rho)$ is an $A$-module $V=V_\ev\oplus V_\od$ together with a $A$-linear morphism of restricted Lie superalgebras $\phi:L\rightarrow \End(V)$ satisfying:
\begin{align}\label{eq:RLRsupermodule-1}
    \phi(x)(av)&=(-1)^{|x||a|}a\phi(x)(v)+\rho(x)(a)v,&~\forall x\in L,~\forall a\in A,~\forall v\in V,\\\label{eq:RLRsupermodule-2}
    \phi(ax)^{p-1}(av)&=a^p\phi(x)^{p-1}(v)+\rho(ax)^{p-1}(a)v,&~\forall x\in L_\ev,~\forall a\in A_\ev,~\forall v\in V,\\\label{eq:RLRsupermodule-3}
    \phi(ax)^{2p-1}(av)&=a^{2p}\phi(x)^{2p-1}(v)+\rho(ax)^{2p-1}(a)v
    \\\nonumber&+\sum_{i=0}^{p-1}\lambda_{i}\rho(ax)^i(a)\rho(ax)^{2p-2-i}(a)\phi(x)(v),&\forall a\in A_\ev,~\forall x\in L_\od~\forall v\in V,
    \end{align}
    where the coefficients $\lambda_{i}$ are defined in Eq. \eqref{eq:lambda}.
Such a pair $(\phi,V)$ is called a \textit{restricted Lie-Rinehart module}.

\subsubsection{Semi-direct product} In this Section, we investigate the semi-direct product of a restricted Lie-Rinehart superalgebra with a representation. We begin with a few technical Lemmas. Recall that if $\bigl(L,[-,-]\bigl)$ is a (ordinary) Lie superalgebra and $(\phi,V)$ a (ordinary) representation, the space $L\oplus V$ is also a Lie superalgebra with the bracket 
\begin{equation}\label{eq:semidirectbracket}
    \bigl[(x+v),(y+w)\bigl]_{\rtimes}:=[x,y]+\phi(x)(w)-(-1)^{|y||v|}\phi(y)(v),~~\forall ~x,y\in L,~\forall~ v,w\in V.
\end{equation}
The Lie superalgebra $(L\oplus V,[-,-]_{\rtimes})$ is denoted by $L\rtimes V$.

\sssbegin{Lemma}\label{lem:semidirectlem1}
    Let $\left(L,[-,-]\right)$ be a Lie superalgebra and $(\phi,V)$ be a representation of $L$. Then for all $n\in\mathbb{N}$, for all $x\in L_\ev$ and all $y\in L$,  we have
    \begin{equation}     \phi\bigl(\ad^n_x(y)\bigl)=\sum_{k=0}^{n}(-1)^k\binom{n}{k}\phi(x)^{n-k}\phi(y)\phi(x)^k.
    \end{equation}
\end{Lemma}

\begin{proof}
    Let $x\in L_\ev$ and $y\in L$. We prove the Lemma by induction on $n\in \mathbb{N}$. The formula is true for $n=0,1$. Then, the inductive step is as follows.
    \begin{align*}
        \phi\bigl(\ad^{n+1}_x(y)\bigl)&=\phi(x)\phi\bigl(\ad_x^n(y)\bigl)-\phi\bigl(\ad_x^n(y)\bigl)\phi(x)\\
        &=\phi(x)\left(\sum_{k=0}^{n}(-1)^k\binom{n}{k}\phi(x)^{n-k}\phi(y)\phi(x)^k\right)-\left(\sum_{k=0}^{n}(-1)^k\binom{n}{k}\phi(x)^{n-k}\phi(y)\phi(x)^k\right)\phi(x)\\
        &=\sum_{k=0}^{n}(-1)^k\binom{n}{k}\phi(x)^{n-k+1}\phi(y)\phi(x)^k-\sum_{k=0}^{n}(-1)^k\binom{n}{k}\phi(x)^{n-k}\phi(y)\phi(x)^{k+1}\\
        &=\phi(x)^{n+1}\phi(y)+\sum_{k=1}^{n}(-1)^k\binom{n}{k}\phi(x)^{n-k+1}\phi(y)\phi(x)^k\\
        &\quad-\sum_{k=0}^{n-1}(-1)^k\binom{n}{k}\phi(x)^{n-k}\phi(y)\phi(x)^{k+1}-(-1)^n\phi(y)\phi(x)^{n+1}\\
        &=\phi(x)^{n+1}\phi(y)+(-1)^{n+1}\phi(y)\phi(x)^{n+1}+\sum_{k=1}^{n}(-1)^k\left(\binom{n}{k}+\binom{n}{k-1}\right)\phi(x)^{n-k+1}\phi(y)\phi(x)^k\\
        &=\phi(x)^{n+1}\phi(y)+(-1)^{n+1}\phi(y)\phi(x)^{n+1}+\sum_{k=1}^{n}(-1)^k\binom{n+1}{k}\phi(x)^{n-k+1}\phi(y)\phi(x)^k\\
        &=\sum_{k=0}^{n+1}(-1)^k\binom{n+1}{k}\phi(x)^{n-k+1}\phi(y)\phi(x)^k,
    \end{align*}
    which finishes the proof.
\end{proof}

\sssbegin{Lemma}\label{lem:semidirectlem2}
    Let $\left(L,[-,-],(-)^{[p]}\right)$ be a restricted Lie algebra, $(\phi,V)$ a restricted representation of $L$ and $L\oplus V$ equipped with the bracket \eqref{eq:semidirectbracket}. Then for all $n\in\mathbb{N}$, all $x+v\in (L\rtimes V)_\ev$ and all $y+w\in L\rtimes V$, we have
    \begin{equation}\label{eq:semidirectlem2n}
        \widetilde{\ad}^n_{(x+v)}(y+w)=\ad_x^n(y)+ \phi(x)^n(w)+\sum_{k=1}^{n}(-1)^k\binom{n}{k}\phi(x)^{n-k}\phi(y)\phi(x)^{k-1}(v),
    \end{equation}
where $\widetilde{\ad}$ refers to the bracket \eqref{eq:semidirectbracket}. In particular, if $n=p$, we have

    \begin{equation}\label{eq:semidirectlem2p}
        \widetilde{\ad}^p_{(x+v)}(y+w)=\ad_x^p(y)+\phi(x)^p(w)-\phi(y)\phi(x)^{p-1}(v).
    \end{equation}

\end{Lemma}

\begin{proof}
    We prove \eqref{eq:semidirectlem2n} by induction on $n\in\mathbb{N}$. We only compute the component in $V$ of the bracket \eqref{eq:semidirectbracket}. Let $x+v\in (L\rtimes V)_\ev$ and $y+w\in L\rtimes V$. Using Lemma \ref{lem:semidirectlem1}, we have
    \begin{align*}
        &\phi(x)\left(\ad^{n}_{(x+v)}(y+w) \right)-\phi\left(\ad_x^n(y)\right)(v)\\
        &=\phi(x)^{n+1}+\sum_{k=1}^{n}(-1)^k\binom{n}{k}\phi(x)^{n-k+1}\phi(y)\phi_{k-1}(v)-\sum_{k=0}^{n}(- 1)^k\binom{n}{k}\phi(x)^{n-k}\phi(y)\phi^{k}(v)\\   
        &=\phi(x)^{n+1}+\sum_{k=1}^{n}\left((-1)^k\binom{n}{k}-(-1)^{k-1}\binom{n}{k-1}\right)\phi(x)^{n-k+1}\phi(y)\phi(x)^{k-1}(v)+(-1)^{n+1}\phi(y)\phi(x)^n(v)\\
        &=\phi(x)^{n+1}+\sum_{k=1}^{n+1}(-1)^{k}\binom{n+1}{k}\phi(x)^{n-k+1}\phi(y)\phi(x)^{k-1}(v),
    \end{align*}
    which finishes the proof.
 \end{proof} \black

\sssbegin{Proposition}\label{prop:rest-semi-direct-Lie} Let $\bigl(L,[-,-],(-)^{[p|2p]}\bigl)$ be a restricted Lie superalgebra and let $(\phi,V)$ be a restricted representation. Then, $L\rtimes V$ is a restricted Lie superalgebra with a $p|2p$-map $(-)^{[p|2p]_{\rtimes}}:L\rtimes V\rightarrow(L\rtimes V)_{\ev}$ satisfying
\begin{equation}
(e_i+v_j)^{[p]_{\rtimes}}=e_i^{[p]}+\phi(e_i)^{p-1}(v_j),
\end{equation} where $(e_i)_i$ forms a basis of $L_\ev$ and $(v_j)_j$ a basis of $V_\ev$.
\end{Proposition}
\begin{proof}
    Let $(e_i)_i$ be a basis of $L_\ev$ and $(v_j)_j$ a basis of $V_\ev$. We define
    $$e_i^{[p]_{\rtimes}}:=e_i^{[p]},\quad\text{and}\quad v_j^{[p]_{\rtimes}}:=0. $$
    Then, we have
    \begin{equation}(e_i+v_j)^{[p]_{\rtimes}}=e_i^{[p]_{\rtimes}}+\sum_{\ell=1}^{p-1}\widetilde{s}_{\ell}(e_i,v_j),\end{equation} where $\widetilde{s_\ell}$ refers to the coefficients defined in Eq. \eqref{eq:def-si} with respect to the bracket \eqref{eq:semidirectbracket}. Using Eq. \eqref{eq:exp-si}, we have 
    \begin{align}
    \sum_{\ell=1}^{p-1}\widetilde{s_\ell}(e_i,v_j)&=\sum_{\underset{x_{p-1}=v_j,~x_p=e_i}{x_k\in\{e_i,v_j\}}}\frac{1}{\sharp(e_i)}[x_1,[x_2,[...,[x_{p-1},x_p]...]_\rtimes]_\rtimes]_\rtimes,
\end{align} where $\sharp(e_i)=\card\{k,~x_k=e_i\}$. All the terms in the above sum are zero except when $\sharp(e_i)=p-1$, therefore we have  
 \begin{equation}\sum_{\ell=1}^{p-1}\widetilde{s_\ell}(e_i,v_j)=-\frac{1}{p-1}\widetilde{\ad}_{e_i}^{p-1}(v_j)=\widetilde{\ad}_{e_i}^{p-1}(v_j)=\phi(e_i)^{p-1}(v_j). \end{equation}
    It remains to show that the map $(-)^{[p]_{\rtimes}}$ satisfies the conditions of Jacobson's Theorem \ref{Jac}, that is, we have $\widetilde{\ad}^p_{e_i+v_j}(y+w)=\widetilde{\ad}_{(e_i+v_j)^{[p]_{\rtimes}}}(y+w),$ for all $y+w\in L\oplus V$ homogeneous. Let $e_i,v_j$ as above and let $y+w$ an homogeneous element of $L\oplus V$. Using Lemma \ref{lem:semidirectlem2}, we have
    \begin{align*}
        &\ad^p_{(e_i+v_j)}(y+w)-\left[(e_i+v_j)^{[p]_{\rtimes}},(y+w)\right]_{\rtimes}\\
        =&\ad_{e_i}^p(y)+\phi(e_i)^p(w)-\phi(y)\phi(e_i)^{p-1}(v_j)-\bigl[e_i^{[p]}+\phi(e_i)^{p-1}(v_j),(y+w)\bigl]_{\rtimes}\\
        =& \ad_{e_i}^{p}(y)+\phi(e_i)^p(w)-\phi(y)\phi(e_i)^{p-1}(v_j)-\bigl(\ad_{e_i^{[p]}}(y)+\phi(e_i^{[p]})(w)-\phi(y)\phi(e_i)^{p-1}(v_j)\bigl)=0.
    \end{align*}
    Therefore, we can apply Jacobson's Theorem \ref{Jac} and $L\rtimes V$ is a restricted Lie superalgebra.
    \end{proof}

\sssbegin{Theorem}\label{thm:semi-direct}
    Let $(A,L,\rho)$ be a restricted Lie superalgebra and let $(\phi,V)$ a representation as in Definition \ref{def:rest-rep}. Suppose that the center of the restricted Lie superalgebra $L\rtimes V$ constructed in Proposition \ref{prop:rest-semi-direct-Lie} is trivial. Then, $(A,L\rtimes V,\widetilde{\rho})$ is a restricted Lie-Rinehart superalgebra, with $\widetilde{\rho}(x+v)=\rho(x),~\forall~x+v\in L\rtimes V$.
\end{Theorem}
\begin{proof}
     Since the center of $L\rtimes V$ is trivial, it follows from Remark \ref{rmk:centerless} and Proposition \ref{prop:rest-semi-direct-Lie} that the $p$-map is given by
     $$(x+v)^{[p]_{\rtimes}}=x^{[p]}+\phi(x)^{p-1}(v),\quad\forall x+v\in L\rtimes V.$$ 
     Let $a\in A$ and $x+v\in L\rtimes V$. We divide the proof into four different cases, depending on the parity of $a$ and $x+v.$

     \underline{The case where $a$ and $x+v$ are even.} We have
      \begin{align*}
                &(ax+av)^{[p]_{\rtimes}}-a^p(x+v)^{[p]_{\rtimes}}-\widetilde{\rho}\bigl(a(x+v)\bigl)^{p-1}(a)(x+v)\\
              =~&\left((ax)^{[p]}+\phi(ax)^{p-1}(av)\right)-a^p\left(x^{[p]}+\phi(x)^{p-1}(v)\right)-\rho(ax)^{p-1}(a)(x+v)\\
              =~&\phi(ax)^{p-1}(av)-a^p\phi(x)^{p-1}(v)-\rho(ax)^{p-1}(a)v=0,
            \end{align*} and thus Eq. \eqref{eq:superhochschild-def-even/even} is satisfied.

     \underline{The case where $a$ is odd and $x+v$ is even.}  Note that in this case, we have
      $$\phi(ax)^2(av)=a\rho(x)(a)^2v+a^2\rho(x)^2(a)v=0.$$ Recall that for $y+w\in (L\oplus V)_{\od}$, we have
     $$ (y+w)^2=\frac{1}{2}\bigl([y,y]_\rtimes+2\phi(y)(w)\bigl)=y^2+\phi(y)(w).$$ Therefore, we have
      $$(y+w)^{[2p]_\rtimes}=y^{[2p]}+\phi(y^2)^{p-1}\circ\phi(y)(w)=y^{[2p]}+\phi(y)^{2p-1}(w).$$ Using the latter identity, we obtain
      $$(ax+av)^{[2p]_\rtimes}=(ax)^{[2p]_\rtimes}+\phi(ax)^{2p-1}(av)=0, $$ and thus Eq. \eqref{eq:superhochschild-def-odd/even} is satisfied.

       \underline{The case where $a$ is even and $x+v$ is odd.} We have
       \begin{align*}
        (ax+av)^{[2p]_{\rtimes}}=&~(ax)^{[2p]}+\phi(ax)^{2p-1}(av)\\
                      =&~a^{2p}x^{[2p]}+\rho(ax)^{2p-1}(a)x+\sum_{i=0}^{p-1}\lambda_{i}\rho(ax)^i(a)\rho(ax)^{2p-2-i}(a)x^2+\phi(ax)^{2p-1}(av)\\
                      =&~a^{2p}\bigl(x^{[2p]}+\phi(x)^{2p-1}(v)\bigl)+\rho(ax)^{2p-1}(x+v)\\
                      &+\sum_{i=0}^{p-1}\lambda_{i}\rho(ax)^i(a)\rho(ax)^{2p-2-i}(a)(x^2+\phi(x)(v))\\
                      =&~a^{2p}(x+v)^{[2p]_{\rtimes}}+\widetilde{\rho}(ax+av)^{2p-1}(a)(x+v)\\
                      &+\sum_{i=0}^{p-1}\lambda_{i}\widetilde{\rho}(ax+av)^i(a)\widetilde{\rho}(ax+av)^{2p-2-i}(a)(x+v)^2,
       \end{align*} and thus Eq. \eqref{eq:superhochschild-def-even/odd} is satisfied.

         \underline{The case where $a$ and $x+v$ are odd.} First, let us show that $\phi(ax)^k(av)=a\rho(x)(a)^kv$, for  all $k\geq 0$. Indeed by induction, we have
        \begin{align*}
            \phi(ax)^{k+1}(av)&=\phi(ax)\bigl(a\rho(x)(a)^{k}v\bigl)\\
            &=a\Bigl(-a\rho(x)(a)^k\phi(x)v+\rho(x)\bigl(a\rho(x)(a)^k)\bigl)v\Bigl)\\
            &=a\rho(x)\bigl(a\rho(x)(a)^k\bigl)v\\
            &=a\rho(x)(a)^{k+1}v-a^2\rho(x)\bigl(\rho(x)(a)^k\bigl)v=a\rho(x)(a)^{k+1}v.
        \end{align*} In particular, we have $\phi(ax)^{p-1}(av)=a\rho(x)(a)^{p-1}v$. It follows that 
       \begin{align*}
       (ax+av)^{[p]_\rtimes}&=(ax)^{[p]_\rtimes}+\phi(ax)^{p-1}(av)\\&=a\bigl(\rho(x)(a)\bigl)^{p-1}x+a\rho(x)(a)^{p-1}v\\&=a\bigl(\widetilde{\rho}(x+v)(a)\bigl)^{p-1}(x+v),
       \end{align*} and thus Eq. \eqref{eq:superhochschild-def-odd/odd} is satisfied. \black
\end{proof}

\section{The restricted enveloping algebra of a restricted Lie-Rinehart superalgebra}\label{sec:UEA}
 The purpose of this section is to construct the restricted universal enveloping algebra of a restricted Lie-Rinehart superalgebra. The non-super case has been investigated in \cite{Do}.

 In this section, all associative superalgebras are assumed to have an \textit{unit}.
\subsection{The restricted enveloping algebra of a restricted Lie superalgebra} Let $L$ be a Lie superalgebra, and let $U(L)$ be its  universal enveloping algebra. In the case where $L$ is restricted, the \textit{restricted} universal enveloping algebra of $L$ is defined (see \cite{Pe}) by
\begin{equation*}
    U_p(L)=U(L)/\langle x^p-x^{[p]},~x\in L_\ev\rangle.
\end{equation*}
Following \cite{Pe}, a Poincaré-Birkhoff-Witt basis of $U_p(L)$ is given by monomials of the form
\begin{equation*}
    e_{i_1}^{k_1}\cdots e_{i_m}^{k_m}f_{j_1}\cdots f_{j_s},\quad i_1<\cdots <i_m,\quad j_1<\cdots<j_s,\quad k_s\leq p-1,
\end{equation*} where $\{e_i\}_i$ is a basis of the even part $L_\ev$ and $\{f_j\}_j$ is a basis of the odd part $L_\od$.
\subsection{The universal enveloping algebra of a Lie-Rinehart superalgebra} We recall the construction of the universal enveloping algebra of a Lie-Rinehart superalgebra, see \cite[Section 5]{La}. Let $(A,L,\rho)$ be a Lie-Rinehart superalgebra. On the superspace $A\oplus L$, a Lie superalgebra structure can be defined by the bracket
\begin{equation}
    [a+x,b+y]=[x,y]+\rho(x)(b)-(-1)^{|a||y|}\rho(y)(a),\quad\forall a+x,b+y\in A\oplus L.
\end{equation} We note the resulting Lie superalgebra by $A\rtimes L$. Note that $A\rtimes L$ is an $A$-module with the action $a(b+x)=ab+ax,~\forall a,b\in A,~\forall x\in L$.

Consider the universal enveloping algebra $U(A\rtimes L)$ of $A\rtimes L$ and denote by $U^+(A\rtimes L)$ the subspace of $U(A\rtimes L)$ spanned by $i(A\rtimes L)$ where $i:A\rtimes L\rightarrow U(A\rtimes L)$ is the (even) injection. Then, the universal enveloping algebra of the Lie-Rinehart superalgebra $(A,L,\rho)$ is defined by
\begin{equation}\label{eq:ideal-J}
    U(A,L)=U^+(A\rtimes L)/J,
    \end{equation} where $J$ is the ideal of $U^+(A\rtimes L)$ spanned by elements $i(a)i(b+x)-i(a(b+x)),~a,b\in A,~x\in L.$ Denote by $\iota_A:A\rightarrow U(A,L)$ and $\iota_L:L\rightarrow U(A,L)$ the natural maps. Then, it is proved in \cite[Section 5]{La} that $\iota_A$ is an injective morphism of associative superalgebras, $\iota_L$ is a morphism of Lie superalgebras, and for all $a\in A$ and all $x\in L$, we have 
\begin{align}
    \iota_A(a)\iota_L(x)&=\iota_L(ax);\\
    \iota_A\bigl(\rho(x)(a)\bigl)&=\iota_L(x)\iota_A(a)-(-1)^{|a||x|}\iota_A(a)\iota_L(x).
\end{align}

The triple $\bigl(U(A,L),\iota_A,\iota_L\bigl)$ satisfies a the following universal property (see \cite[Section 5]{La}). Let $B$ be an associative superalgebra, let $j_A:A\rightarrow B$ be an even morphism of associative superalgebras and let $j_L:L\rightarrow B$ be an even morphism of Lie superalgebras (where $B$ is seen as a Lie superalgebra under the commutator) satisfying for all $a\in A$ and all $x\in L$ the conditions
\begin{equation}
    j_L(ax)=j_A(a)j_L(x);\quad\text{and}\quad j_A\bigl(\rho(x)(a) \bigl)=j_L(x)j_A(a)-(-1)^{|a||x|}j_A(a)j_L(x).
\end{equation}
Then, there exists an unique morphism of associative superalgebras $\phi:U(A,L)\rightarrow B$ such that $\phi\circ \iota_L=j_L$ and $\phi\circ \iota_A=j_A$.

\subsection{The universal enveloping algebra of a restricted Lie-Rinehart superalgebra} We follow the construction described by Dokas in \cite[Section 2]{Do} minding the sign rule. Let $(A,L,\rho)$ be a \textit{restricted} Lie-Rinehart superalgebra. Let $\pi_1:U^+(A\rtimes L)\to U(A,L)=U^+(A\rtimes  L)/J$ be the projection. We define the restricted universal enveloping algebra of $(A,L,\rho)$ by $$U_p(A,L):=U(A,L)/\langle \pi_1\circ i(x^{[p]})-(\pi_1\circ i(x))^p,~x\in L_\ev\rangle.$$ 

We denote by $\pi_2:U(A,L)\to U_p(A,L)$ the projection. Moreover, let $i_A=\pi_2\circ \iota_A$ and $i_L=\pi_2\circ \iota_L$.  For all $a\in A$ and all $x\in L$, we have 
\begin{align}
    i_A(a)i_L(x)&=i_L(ax);\\
    i_A\bigl(\rho(x)(a)\bigl)&=i_L(x)i_A(a)-(-1)^{|a||x|}i_A(a)i_L(x).
\end{align}

We have the following universal property.
\sssbegin{Proposition}\label{prop:univ-prop}
     Let $B$ be an associative superalgebra, let $j_A:A\rightarrow B$ be an even morphism of associative superalgebras and let $j_L:L\rightarrow B$ be an even morphism of restricted Lie superalgebras \textup{(}where $B$ is seen as a restricted Lie superalgebra under the commutator and the $p$-th power on even elements\textup{)} satisfying for all $a\in A$ and all $x\in L$ the conditions
\begin{equation}\label{eq:univ-j-rest}
    j_L(ax)=j_A(a)j_L(x);\quad\text{and}\quad j_A\bigl(\rho(x)(a) \bigl)=j_L(x)j_A(a)-(-1)^{|a||x|}j_A(a)j_L(x).
\end{equation}
Then, there exists a unique morphism of associative superalgebras $\psi:U_p(A,L)\rightarrow B$ such that $\psi\circ i_L=j_L$ and $\psi\circ i_A=j_A$.
\end{Proposition}
\begin{proof}
       Consider a triple $(B,j_A,j_L)$ satisfying \eqref{eq:univ-j-rest}. We define a linear map $$\varphi:A\rtimes L \to B,\quad (a+x)\longmapsto j_A(a)+j_L(x).$$   
       The map $\varphi$ is a morphism of Lie superalgebras. Indeed, for all $a,b\in A$ and all $x,y\in L$, we have
       \begin{align*}
        &\varphi(a+x)\varphi(b+y)-(-1)^{|a||y|}\varphi(b+y)\varphi(a+x)\\
        =&~(j_A(a)+j_L(x))(j_A(b)+j_L(y))-(-1)^{|a||y|}(j_A(b)+j_L(y))(j_A(a)+j_L(x))\\
        =&~j_A(a)j_L(y)+j_L(x)j_A(b)+j_L(x)j_L(y)-(-1)^{|a||y|}\bigl(j_A(b)j_L(x)+j_L(y)j_A(a)+j_L(y)j_L(x)\bigl)\\
        =&~j_A\bigl(\rho(x)(b)\bigl)-(-1)^{|a||y|}j_A\bigl(\rho(y)(a)\bigl)+j_L\bigl([x,y]\bigl)\\
        =&~\varphi\bigl([a+x,b+y]\bigl).
       \end{align*}
By the universal property of $U(A\rtimes L)$, there exists a unique associative superalgebras morphism $\widetilde{\varphi}:U^+(A\rtimes L)\to B$ such that $\widetilde{\varphi}\circ i=\varphi.$
       Moreover, for all $a,b\in A$ and all $x\in L$, we have
       \begin{align*}
            \widetilde{\varphi}\bigl(i(a)i(b+x)-i(a(b+x))\bigl)&=\widetilde{\varphi}( i(a))\widetilde{\varphi}(i(b+x))-\widetilde{\varphi}(i(ab+ax))\\
            &=\varphi(a)\varphi(b+x)-\varphi(ab+ax)\\
            &=j_A(a)\bigl(j_A(b)+j_L(x)\bigl)-j_A(ab)-j_L(ax)=0.
       \end{align*} It follows that the ideal $J$ (see Eq. \eqref{eq:ideal-J}) is contained in the kernel of the map $\widetilde{\varphi}$. Thus, $\widetilde{\varphi}$ induces a unique associative superalgebras morphism $\widetilde{\psi}:U(A,L)\to B$. 
Since $j_L$ is a morphism of restricted Lie superalgebras, we have
\begin{align*}
\widetilde{\psi}\bigl(\pi_1\circ i(x)^p\bigl)-\widetilde{\psi}\bigl(\pi_1\circ i(x^{[p]})\bigl)&=\bigl(\widetilde{\psi}\circ\pi_1\circ i(x)\bigl)^p-\widetilde{\psi}\circ\pi_1\circ i(x^{[p]})\\
&=\varphi(x)^p-\varphi(x^{[p]})\\
&=j_L(x)^p-j_L(x^{[p]})=0.
\end{align*}
Thus, $\widetilde{\psi}$ induces a unique associative superalgebras morphism $\psi:U_p(A,L)\rightarrow B$. We have the following commutative diagram:
       \begin{equation*}
\begin{tikzcd}
A\oplus L \arrow[rrrr, "\varphi"] \arrow[rrd, "i"', hook] &  &                                                                         &  & B                                       \\
 &  & U^+(A\rtimes L) \arrow[rru, "\widetilde{\varphi}"] \arrow[dd, "\pi_1"']  &  &                                         \\
   &  &                                                                         &  &                                         \\
    &  & {U(A,L)} \arrow[rruuu, "\widetilde{\psi}"', dotted] \arrow[rr, "\pi_2"] &  & {U_p(A,L)} \arrow[uuu, "\psi"', dashed]
\end{tikzcd}
\end{equation*} 
Moreover, the map $\psi$ satisfies $\psi\circ i_A(a)=\varphi(a)=j_A(a)~\forall a\in A$ as well as $\psi\circ i_L(x)=\varphi(x)=j_L(x)~\forall x\in L.$ 
\end{proof}


\appendix
\renewcommand{\thesection}{\Alph{section}}
\renewcommand{\thesubsection}{\thesection.\arabic{subsection}}
\section{Appendix: proof of Proposition \ref{prop:p-divides-Gamma2pp}}\label{sec:gamma2pp}

Here we detail the proof of Proposition \ref{prop:p-divides-Gamma2pp}.
\subsection{Preliminary results} Let us start with some preliminary results.
\sssbegin{Lemma}
   Let $k\geq 3$ and $2\leq j\leq k-1$. We have
   \begin{equation}\label{eq:UUU}
    \Gamma_{k,j}=x_0^2\Gamma_{k-2,j-1}+x_0^2\delta^2(\Gamma_{k-2,j})+x_0x_1\bigl(\delta(\Gamma_{k-2,j})-(-1)^{k-j}\Gamma_{k-2,j-1}\bigl).
   \end{equation}
\end{Lemma}
\begin{proof}Let $k\geq 3$ and $0\leq j\leq k-1$. It follows from the recursive definition of the $\Gamma$'s (see Lemma \ref{lem:recursive-Gamma}) that $|\Gamma_{k,j}|=k-j\mod 2$. Thus, we have
\begin{align}\label{eq:1}
    \Gamma_{k,j}&=(-1)^{k-j}x_0\Gamma_{k-1,j-1}+x_0\delta(\Gamma_{k-1,j}),\\
\label{eq:2}
    \Gamma_{k-1,j-1}&=(-1)^{k-j}x_0\Gamma_{k-2,j-2}+x_0\delta(\Gamma_{k-2,j-1}),
    \end{align} and
    \begin{equation}\label{eq:3}
    \Gamma_{k-1,j}=(-1)^{k-j+1}x_0\Gamma_{k-2,j-1}+x_0\delta(\Gamma_{k-2,j}).
\end{equation}
By multiplying Eq. \eqref{eq:2} by $(-1)^{k-j}x_0$ and substituting into Eq. \eqref{eq:1}, we have 
\begin{equation}\label{eq:F}
   \Gamma_{k,j}=x_0^2\Gamma_{k-2,j-2}+(-1)^{k-j}x_0^2\delta\bigl(\Gamma_{k-2,j-1}\bigl)+x_0\delta\bigl(\Gamma_{k-1,j}\bigl).
    \end{equation} 
Applying the derivation $\delta$ to \eqref{eq:3} yields
    \begin{align}\label{eq:4}
        \delta(\Gamma_{k-1,j})&=(-1)^{k-j+1}\delta(x_0\Gamma_{k-2,j-1})+\delta\bigl(x_0\delta(\Gamma_{k-2,j})\bigl)\\\nonumber
        &=(-1)^{k-j+1}\bigl(x_1\Gamma_{k-2,j-1}+x_0\delta(\Gamma_{k-2,j-1})\bigl)+x_1\delta(\Gamma_{k-2,j})+x_0\delta^2(\Gamma_{k-2,j})\bigl).
    \end{align}
    By substituting \eqref{eq:4} into \eqref{eq:F}, we finally obtain \eqref{eq:UUU}.
\end{proof}

\sssbegin{Remark}\label{rmk:P-zero} We may extend the definition of the $\Gamma_{k,j}$ to the half-plane $\{(k,j),~k\geq 0,~j\in\Z\}$ by setting $\Gamma_{0,0}=1$, $\Gamma_{0,j}=0~(j\neq 0)$ and
$$\Gamma_{k,j}=x_0\delta(\Gamma_{k-1,j})+(-1)^{k-j}x_0\Gamma_{k-1,j-1},~\forall k\geq 1,~\forall j\in \Z.$$
Note that by induction on $k$, for every $k \geq 1,$ we have $\Gamma_{k,j} = 0$ whenever $j \leq 0$ or $j \geq k + 1.$
\end{Remark}

\subsubsection{Notations} Denote by $r = k - j$ the quantity we will refer to as the \textit{diagonal}. 
It follows from the recursion (see Lemma~\ref{lem:recursive-Gamma}) that each monomial appearing in $\Gamma_{k,j}$ has the form
\[
x_0^{k-s} x_{a_1} \cdots x_{a_s},
\]
where $1 \leq a_1 \leq \cdots \leq a_s$ and $a_1 + a_2 + \cdots + a_s = k - j$.

Let $\mu = (a_1, a_2, \ldots, a_s)$ represents such a sequence, which we will call a \textit{shape}. We denote its \textit{length} by $|\mu|=s$. In the sequel, we shall use the notation $x^\mu=x_{a_1} \cdots x_{a_s}$. For $t\geq 0$, define the set of $t$\textit{-shapes} by
\begin{align*}S^{t}&:=\bigl\{\mu=(a_1\le\cdots\le a_{|\mu|}),\ a_1+\cdots+a_{|\mu|}=t\bigl\}
,\quad\text{and let}\\
X^t&:=\bigl\{(\mu,u),~\mu\in S^t,~1\leq u\leq|\mu|\bigl\}. \end{align*}
Let $\mu=(a_1,a_2,\cdots,a_{|\mu|})$ be a finite sequence. We denote by $\sort(\mu)$ its nondecreasing reordering and by $\mu\cup \{k\}$ the sequence $(k,a_1,a_2,\cdots,a_{|\mu|})$. For $u\in\{1,\cdots,|\mu|\}$ and $i\geq 1$, we define
$$\inc_u^{(i)}\mu:=(a_1,\cdots,a_{u-1},a_u+i,a_{u+1},\cdots,a_{|\mu|})$$ the sequence  obtained by increasing the $u$-th entry of $\mu$ by $i$. We define the following maps.
$$\begin{array}{cccccl}
    \Phi_i^t&:&S^t\rightarrow S^{t+i},& \Phi_i^t(\mu)&=&\sort(\mu\cup\{i\});\\
    (\Phi_i^t)^{-1}&:&S^{t+i}\rightarrow \mathcal{P}(S^t),&(\Phi_i^t)^{-1}(\nu)&=&\bigl\{\mu\in S^t,~\Phi_i^t(\mu)=\nu\bigl\}.\\
    \Psi_i^t&:&X^t\rightarrow S^{t+i},& \Psi_i^t(\mu,u)&=&\sort(\inc_u^{(i)}\mu);\\
    (\Psi_i^t)^{-1}&:&S^{t+i}\rightarrow \mathcal{P}(X^t),&(\Psi_i^t)^{-1}(\nu)&=&\bigl\{(\mu,u)\in X^t,~\Psi_i^t(\mu,u)=\nu\bigl\}.  
\end{array}$$
We also define restriction as follows, for $U\subset S^t$.
$$\begin{array}{cccccl}
   
    (\Phi_i^t|_U)^{-1}&:&S^{t+i}\rightarrow \mathcal{P}(U),&(\Phi_i^t|_U)^{-1}(\nu)&=&(\Phi_i^t)^{-1}(\nu)\cap U.\\
   
    (\Psi_i^t|_U)^{-1}&:&S^{t+i}\rightarrow \mathcal{P}(X^t_U),&(\Psi_i^t|_U)^{-1}(\nu)&=&\bigl\{ (\mu,u)\in X_U^t,~\Psi_i^t(\mu,u)=\nu\bigl\}, 
\end{array}$$ where $X^t_U=\bigl\{(\mu,u)\in X^t,~\mu\in U\bigl\}.$

For all $t\geq 0$, we define the \textit{Leibniz sign} for a pair $(\mu,u)\in X^t$ by
$$\varepsilon(\mu,u)=(-1)^{\sharp\{i<u,~a_i~\text{odd}\}},$$
where $\mu= (a_1, a_2, \ldots, a_s).$

Let $t\geq 0$, $U\subset S^t$ and $N_{\alpha,\beta}(m):=\alpha m+\beta,~\alpha,\beta\in\Z$. A \textit{partial homogeneous sum} is 
\begin{equation}\label{eq:partial-sum}H_{t, N}^{[U]}(m) = \sum_{\mu \in U} C_{\mu}(m) \cdot  x_0^{N_{\alpha,\beta}(m) - |\mu|} x^\mu,\end{equation} where each $C_\mu(m)$ is an integer valued function in $m$ (it is supposed to be zero whenever $N_{\alpha,\beta}(m)-|\mu|<0$). 
Every monomial in $H_{t, N}^{[U]}(m)$ is of total degree $N_{\alpha,\beta}(m).$ 

Finally, if $P$ is a polynomial and $X$ is a monomial that appears in $P$, we denote by $[X]P$ the coefficient of $X$ in $P$.

\subsubsection{Some identities}
 Using the previously defined notations, for any shape $\mu$, we have 
 \begin{equation*}
    \delta(x^\mu)=\sum_{u=1}^{|\mu|}\varepsilon(\mu,u)x^{\inc_u^{(1)}(\mu)}.
 \end{equation*}
Now, let $\mu\in S^t$ and let $A\in\Z$ be a fixed integer. We have
\begin{equation*}
    \delta(x_0^Ax^\mu)=\delta(x_0^A)x^\mu+x_0^A\delta(x^\mu)=Ax_0^{A-1}x_1x^\mu + \sum_{u=1}^{|\mu|} \varepsilon(\mu,u)x^A x^{\inc_u^{(1)}(\mu)}.
\end{equation*}
Consider a partial homogeneous sum $H_{t, N}^{[U]}(m)$ as defined in Eq. \eqref{eq:partial-sum}. We have
$$ \delta\bigl(H_{t, N}^{[U]}(m)\bigl)=\sum_{\mu\in U} C_\mu(m)\Big((N_{\alpha,\beta}(m)-|\mu|)\,x_1\,x_0^{\,N_{\alpha,\beta}(m)-|\mu|-1}\,x^\mu
\ +\ \sum_{u=1}^{|\mu|}\varepsilon(\mu,u)\,x_0^{\,N_{\alpha,\beta}(m)-|\mu|}\,x^{\mathrm{inc}^{(1)}_u(\mu)}\Big)
.$$ Thus, it follows that for $\nu \in S^{t+1}$, we have
\begin{align}\label{eq:D3}
\big[x_0^{\,N_{\alpha,\beta}(m)-|\nu|}x^\nu\big]\ \delta\big(H^{[U]}_{t,N}(m)\big)
=&\sum_{\mu\in(\Phi_1^t|_{U})^{-1}(\nu)}\!\!\big(N_{\alpha,\beta}(m)-|\mu|\big)\,C_\mu(m)
\\\nonumber&+\;\sum_{(\mu,u)\in(\Psi_1^t|_{U})^{-1}(\nu)}\!\!\varepsilon(\mu,u)\,C_\mu(m).
\end{align}
Similarly, we have
\[
\delta^{2}\big(H^{[U]}_{t,N}(m)\big)
=\sum_{\mu\in U} C_\mu(m)\Big((N_{\alpha,\beta}(m)-|\mu|)\,x_2\,x_0^{\,N_{\alpha,\beta}(m)-|\mu|-1}\,x^\mu
\ +\ \sum_{u=1}^{|\mu|}x_0^{\,N_{\alpha,\beta}(m)-|\mu|}\,x^{\mathrm{inc}^{(2)}_u(\mu)}\Big). 
\] It follows that for $\nu \in S^{t+2}$, we have
\[
\big[x_0^{\,N_{\alpha,\beta}(m)-|\nu|}x^\nu\big]\ \delta^{2}\!\big(H^{[U]}_{t,N}(m)\big)
=\sum_{\mu\in(\Phi_2^t|_{U})^{-1}(\nu)}\!\!\big(N_{\alpha,\beta}(m)-|\mu|\big)\,C_\mu(m)
\;+\;\sum_{(\mu,u)\in(\Psi_2^t|_{U})^{-1}(\nu)}\!\! C_\mu(m);
\] and 
\begin{align}\label{eq:D4}
\big[x_0^{\,N_{\alpha,\beta}(m)+2-|\nu|}x^\nu\big]\ x_0^{2}\,\delta^{2}\!\big(H^{[U]}_{t,N}(m)\big)
=&\sum_{\mu\in(\Phi_2^t|_{U})^{-1}(\nu)}\!\!\big(N_{\alpha,\beta}(m)-|\mu|\big)\,C_\mu(m)
\;\\\nonumber&+\;\sum_{(\mu,u)\in(\Psi_2^t|_{U})^{-1}(\nu)}\!\! C_\mu(m).
\end{align}

Let $r\geq 0$ be an integer and consider a $r$-shape $\lambda$.
We define a sequence $P_\lambda(m), \; m\geq 0$ as follows.
$$P_\lambda(m) = \big[x_0^{2m - |\lambda|} x^{\lambda}\big] \Gamma_{2m, 2m-r},\quad m \geq \lceil \frac{r+1}{2} \rceil$$
and 
$$P_\lambda(m) = 0,\quad 0 \leq m < \lceil \frac{r+1}{2}\rceil.$$
Analogously, we define a sequence $Q_\lambda(m)$ as follows.
$$Q_\lambda(m) = \big[x_0^{2m+1-|\lambda|}x^{\lambda}\big] \Gamma_{2m+1, 2m+1-r},\quad m \geq \lceil \frac{r}{2} \rceil,$$
and 
$$Q_\lambda(m) = 0,\quad 0 \leq m < \lceil\frac{r}{2} \rceil.$$

\sssbegin{Lemma}\label{lem:polym}
   Let $r\geq 1$. For any $r$-shape $\lambda$, the sequences $P_\lambda(m)$ and $Q_\lambda(m)$ are polynomials in $m$. 
\end{Lemma}

\begin{proof}
        Recall that if $f=\sum_{i=0}^nb_ix^i$ is a polynomial of degree $n$ in one variable $x$, its finite difference $\Delta f(x)=f(x)-f(x-1)$ is a polynomial of degree $n-1$ and the coefficient of $x^{n-1}$ is given by $nb_n$.
        
Using the coefficient functional defined by
$$\big[x_0^n x^\lambda\big]\bigl(\sum_{(m, \mu)} c_{m, \mu} x_0^mx^\mu\bigr) = c_{n, \lambda} \textbf{1}_{n \geq 0,}$$
we have ($r=|\lambda|$)
$$P_\lambda(m) = [x_0^{2m-|\lambda|}x^\lambda]\Gamma_{2m, 2m-r} \; \forall \; m \geq0$$
and
$$Q_\lambda(m) = [x_0^{2m+1-|\lambda|}x^\lambda]\Gamma_{2m+1, 2m+1-r} \; \forall \; m  \geq 1.$$
        
  We will prove the Lemma by induction on $r$. In the case where $r=1$, we have
    \begin{align*}
        \Gamma_{k,k-1}&=-x_0\Gamma_{k-1,k-2}+x_0\delta(\Gamma_{k-1,k-1})\\
        &=-x_0\Gamma_{k-1,k-2}+(k-1)x_1x_0^{k-1}\\
        &=x_0^2\Gamma_{k-2,k-3}-(k-2)x_1x_1x_0^{k-1}+(k-1)x_1x_0^{k-1}\\
        &=x_0^2\Gamma_{k-2,k-3}+x_1x_0^{k-1}.
    \end{align*}
    Therefore, we have a 2-step recursion and we obtain $ \Gamma_{k,k-1}=mx_0^{k-1}x_1$. The case $r=1$ follows.
Let $r\ge2$.

\underline{The case where $k=2m$}. We have $j=2m-r$.
Extracting the
coefficient of $x_0^{2m-|\lambda|}x^\lambda$ in Eq. \eqref{eq:F},   we obtain for all $m\geq 1$:
\[
\Delta P_\lambda(m):=A_\lambda(m)+B_\lambda(m),
\]where
\begin{align*}
    A_\lambda(m)&=\big[x_0^{\,2m-|\lambda|}x^\lambda\big]\ x_0\delta( \Gamma_{2m-1,2m-r});\\
    B_\lambda(m)&=(-1)^r\big[x_0^{\,2m-|\lambda|}x^\lambda\big]\ x_0^{2}\delta( \Gamma_{2m-2,2m-r-1}).
\end{align*}
The two lower rows in the partial homogeneous sum form are given by
\begin{align}
\Gamma_{2m-1,2m-r}&=\sum_{\mu\in S^{r-1}}Q_\mu(m-1)\;x_0^{\,2m-1-|\mu|}\,x^\mu,
\\\label{eq: pmu-even}
\Gamma_{2m-2,2m-r-1}&=\sum_{\mu\in S^{r-1}}P_\mu(m-1)\;x_0^{\,2m-2-|\mu|}\,x^\mu.
\end{align}
Thus, we have $\Gamma_{2m-1,2m-r}=H^{[S^{r-1}]}_{t=r-1,N_{\alpha,\beta}(m)=2m-1}(m-1)$ with coefficients $Q_\mu(m-1)$, while 
$\Gamma_{2m-2,2m-r-1}=H^{[S^{r-1}]}_{t=r-1,N_{\alpha,\beta(m)}=2m-2}(m-1)$ with coefficients $P_\mu(m-1)$.

By applying Eq. \eqref{eq:D3} to each $x_0\delta(-)$
for $\nu=\lambda\in S^{r}$, we have
\[
\begin{aligned}
A_\lambda(m)
&=\sum_{\mu\in(\Phi_1^{\,r-1})^{-1}(\lambda)}\!\!\big(2m-1-|\mu|\big)\,Q_\mu(m-1)
\;+\;\sum_{(\mu,u)\in(\Psi_1^{\,r-1})^{-1}(\lambda)}\!\!\varepsilon(\mu,u)\,Q_\mu(m-1),\\[2mm]
B_\lambda(m)
&=(-1)^r\Bigg(\sum_{\mu\in(\Phi_1^{\,r-1})^{-1}(\lambda)}\!\!\big(2m-1-|\mu|\big)\,P_\mu(m-1)
\;+\;\sum_{(\mu,u)\in(\Psi_1^{\,r-1})^{-1}(\lambda)}\!\!\varepsilon(\mu,u)\,P_\mu(m-1)\Bigg).
\end{aligned}
\]
By the induction hypothesis, each $P_\mu(-)$ and $Q_\mu(-)$ is a polynomial in $m$; the only additional
$m$-dependence here is the linear factor $(2m-1-|\mu|)$. Since maps $\Phi_1^t, \Psi_1^t, (\Phi_1^t)^{-1}, (\Psi_1^t)^{-1}$ are independent of $m$, expressions $A_\lambda(m)$ and $B_\lambda(m)$ will be polynomials
in $m$, and so is $\Delta P_\lambda(m)$ for all $m\ge1$. 
Define $R_\lambda(m):=\Delta P_\lambda(m)$ and
\[
S_\lambda(m):=\sum_{t=1}^{m}R_\lambda(t)\qquad(m\in\mathbb Z_{\ge0}).
\]
Since $R_\lambda(t)$ is a polynomial in $t$, $S_\lambda(m)$ is a polynomial in $m$. Moreover, we have $S_\lambda(0)=0$.
We also have $P_\lambda(0)=0$ and
$\Delta S_\lambda(m)=R_\lambda(m)=\Delta P_\lambda(m)$ for $m\ge1$. By induction on $m$,
$P_\lambda(m)=S_\lambda(m)$ for all $m\ge0$. It follows that $P_\lambda(m)$ is a polynomial.

 \underline{The case where $k=2m+1$.} We have
\begin{align}\label{eq: pmu-odd}
\Gamma_{2m,2m+1-r}&=\sum_{\mu\in S^{r-1}}P_\mu(m-1)\,x_0^{\,2m-|\mu|}x^\mu,\\
\Gamma_{2m-1,2m+1-r-1}&=\sum_{\mu\in S^{r-1}}Q_\mu(m-1)\,x_0^{\,2m-1-|\mu|}x^\mu.
\end{align}
Similar to the even case, it follows that $\Delta Q_\lambda(m)$ is a polynomial in $m$. Consequently, $Q_\lambda(m)$ is a polynomial in $m$ as $Q_\lambda(0)=0$.
\end{proof}

 A $r$-shape $\lambda$ is called \textit{packed} if $\lambda = (2^{\frac{r}{2}})$ for an even $r$ and $\lambda = (1, 2^{\frac{r-1}{2}})$ for an odd $r.$

\sssbegin{Lemma}\label{lem:A16}
    Let $\lambda$ be a $r-$shape. Then, we have $\deg P_\lambda(m)\leq r$ and  $\deg Q_\lambda(m) \leq r.$ Moreover, if $\deg P_\lambda(m) = r$, 
\end{Lemma}
\begin{proof}
    We will prove the statement by induction on $r$. The induction base is given by cases $r = 1$ and $r = 2$. Recall the the case $r=1$ was computed in the proof of Lemma \ref{lem:polym}  and a direct computation shows that $\Gamma_{2m, 2m-2} = m(m-1)x_0^{2m-1}x_2, \; \Gamma_{2m+1,2m-1} = (m-1)^2 x_0^{2m}x_2.$
    
    First, we prove that $\deg P_\lambda(m), \deg Q_\lambda(m) \leq r.$ Without loss of generality, assume that $k = 2m.$ By comparing the coefficient of $x_0^{2m-|\lambda|}x^\lambda$ in Eq. \eqref{eq:UUU}, we have
\begin{align}\label{eq:above}\Delta P_\lambda(m) =& \big[ x_0^{2m-|\lambda|}x^\lambda\bigl] \; x_0^2 \delta^2 \Gamma_{2m-2, j} + \bigl[x_0^{2m-|\lambda|} x^{\lambda}\bigr]\;x_0x_1\delta(\Gamma_{2m-2,j})\\\nonumber &- (-1)^{r} \bigl[x_0^{2m-|\lambda|}x^{\lambda}\bigr] \Gamma_{2m-2, j-1}.\end{align}
We compute the right-hand side of \eqref{eq:above}. Let $\lambda\in S^r$. It follows from Eqs. \eqref{eq: pmu-even} that
\begin{equation}
    \Gamma_{2m-2, 2m-r} =  \sum_{\mu \in S^{r-2}} P_\mu(m-1) x_0^{2m-2 - |\mu|}x^\mu.
\end{equation}
Using \eqref{eq:D4} (with $t = r-2, N_{\alpha,\beta} = 2m-2, U = S^{r-2}$), we have $$ H_{t,N}^{|U|}(m) = \Gamma_{2m-2,2m-r} = \sum_{\mu \in S^{r-2}} P_\mu(m-1)x_0^{2m-2 - |\mu|}x^\mu$$ (so $C_\mu(m) = P_\mu(m-1)$) and for $\nu = \lambda$, it follows that
\begin{align}\label{eq:F1}\bigl[x_0^{2m-2 + 2 - |\lambda|} x^\lambda\bigr] \bigl(x_0^2\delta^2(\Gamma_{2m-2, 2m-r})\bigr) =& \sum_{\mu \in \big(\Phi_2^{r-2}|_{S^{r-2}}\big)^{-1}(\lambda)} (2m - 2 - |\mu|) P_\mu(m-1) \\\nonumber&+\sum_{(\mu, u) \in \big(\Psi_2^{r-2}|_{S^{r-2}}\big)(\lambda)} P_\mu(m-1).\end{align} 
Similarly, we have
\begin{equation}\label{eq:F2}
\bigl[x_0^{\,2m-|\lambda|}x^\lambda\bigr]
\bigl(x_0x_1\,\delta\,\Gamma_{2m-2, 2m-r}\bigr)
=
\sum_{\nu \in (\Phi_1^{\,r-1}|_{U^{\,r-1}_{\mathrm{no1}}})^{-1}(\lambda)}
\;
\sum_{(\mu,u)\in (\Psi_1^{\,r-2})^{-1}(\nu)}
\varepsilon(\mu,u)\,P_\mu(m-1),
\end{equation} where $U^{\,r-1}_{\mathrm{no1}} := \{\,\nu \in S^{r-1},\; 1 \notin \nu\,\},$ and
\begin{align}\label{eq:F3}
\big[x_0^{2m-|\lambda|}x^\lambda\big]\Big( x_0x_1\Gamma_{2m-2, 2m-1-r}\Big) &= [x_0^{2m-|\lambda|}x^\lambda] \Big(\sum_{\mu \in S^{r-1}} P_\mu(m-1) x_0^{2m-1-|\mu|}x_1x^{\mu}\Big)\\ \nonumber
 &= \sum_{\mu \in \big(\Phi_1^{r-1}|_{\{\mu \in S^{r-1} : 1 \notin \mu\}}\big)^{-1}(\lambda)} P_\mu(m-1).
\end{align}
Using \eqref{eq:F1}, \eqref{eq:F2} and \eqref{eq:F3}, we have
\begin{align}\label{eq:RHS}\nonumber\Delta P_\lambda(m) =& \sum_{\mu \in \big(\Phi_2^{r-2}|_{S^{r-2}}\big)^{-1}(\lambda)} (2m - 2 - |\mu|) P_\mu(m-1) +\sum_{(\mu, u) \in \big(\Psi_2^{r-2}|_{S^{r-2}}\big)(\lambda)} P_\mu(m-1) \;\\&+ \sum_{\nu \in (\Phi_1^{\,r-1}|_{U^{\,r-1}_{\mathrm{no1}}})^{-1}(\lambda)}
\;
\sum_{(\mu,u)\in (\Psi_1^{\,r-2})^{-1}(\nu)}
\varepsilon(\mu,u)\,P_\mu(m-1)\\\nonumber &+ \sum_{\mu \in \big(\Phi_1^{r-1}|_{\{\mu \in S^{r-1} : 1 \notin \mu\}}\big)^{-1}(\lambda)} P_\mu(m-1).
\end{align}
Note that we have $\deg\bigl( (2m-2-|\mu|)P_\mu(m-1)\bigl)\leq r-2 + 1$ since $\mu \in \big(\Phi_2^{r-2}|_{S^{r-2}}\big)^{-1}(\lambda).$ Analogously, every term on the right hand side has degree maximum of $r-1.$ Thus, $\Delta P_\lambda(m)$ has degree $\leq r-1.$ This proves that $\deg P_\lambda(m) \leq r.$ The similar treatment for $Q_\lambda(m)$ proves that $\deg Q_\lambda(m) \leq r.$

Let us show that the maximal degree is obtained for packed $\lambda$'s. 
\begin{equation}\label{eq:prop1a}
\begin{minipage}{0.85\textwidth}
Suppose that $r$ is odd. The packed $(r-2)-$shape is $\lambda_{\text{packed}} = (1,2^{\frac{r-3}{2}}).$ Suppose that for some $\lambda \in S^{r}$, $\lambda_{(r-2)\text{ packed}} \in \Big(\Phi_2^{r-2}|_{S^{r-2}}\Big)^{-1}(\lambda).$ Then $\lambda$ is $r-$packed, i.e. $\lambda = (1, 2^{\frac{r-1}{2}}).$
\end{minipage}
\end{equation}
\noindent\textit{Proof of Statement \eqref{eq:prop1a}.} By definition of $\Phi_2^t(\mu)$ and $\big(\Phi_2^t\big)^{-1}$, we have:
$$\lambda = \Phi_2^t(\lambda_{(r-2) \text{ packed}}) = \text{sort}((1, 2^{\frac{r-3}{2}}) \cup \{2\}) = (1,  2^{\frac{r-1}{2}}) \quad \square.$$
Similarly, we obtain the following statements.
\begin{equation}\label{eq:prop1b}
\begin{minipage}{0.85\textwidth}
Suppose that $r$ is even. The packed $(r-2)-$shape is $\lambda_{\text{packed}} = (2^{\frac{r-2}{2}}).$ Suppose that for some $\lambda \in S^{r}$, $\lambda_{(r-2)\text{ packed}} \in \Big(\Phi_2^{r-2}|_{S^{r-2}}\Big)^{-1}(\lambda).$ Then $\lambda$ is $r-$packed, i.e. $\lambda = (2^{\frac{r}{2}}).$
\end{minipage}
\end{equation}

\begin{equation}\label{eq:prop2a}
\begin{minipage}{0.85\textwidth}
Suppose that $r$ is odd. The packed $(r-1)-$shape is $\lambda_{(r-1) \text{ packed}} = (2^{\frac{r-1}{2}}).$ Suppose that for some $\lambda \in S^{r}, \lambda_{(r-1) \text{ packed}} \in \Big(\Phi_1^{r-1}|_{\{\mu \in S^{r-1: 1\notin \mu}\}}\Big)^{-1}(\lambda).$ Then $\lambda$ is $r-$packed, i.e. $\lambda = (1,2^{\frac{r-1}{2}}).$ 
\end{minipage}
\end{equation}

\begin{equation}\label{eq:prop2b}
\begin{minipage}{0.85\textwidth}
Suppose that $r$ is even. The packed $(r-1)-$shape is $\lambda_{(r-1) \text{ packed}} = (1,2^{\frac{r-2}{2}}).$ Then no $\lambda \in S^{r}$ satisfies $ \lambda_{(r-1) \text{ packed}} \in \Big(\Phi_1^{r-1}|_{\{\mu \in S^{r-1: 1\notin \mu}\}}\Big)^{-1}(\lambda).$ 
\end{minipage}
\end{equation}

We have to prove that if $\lambda \in S^r$ and $\lambda$ is not packed, then $\deg P_\lambda < r.$ 
If $\lambda_{(r-2) \text{-packed}} \in \Big(\Phi_{2}^{r-2}|_{S^{r-2}}\Big)^{-1}(\lambda),$ then it follows from \eqref{eq:prop1a} and \eqref{eq:prop1b} that $\lambda$ is $r-$packed, a contradiction. Thus, $\lambda_{(r-2) \text{-packed}} \notin \Big(\Phi_{2}^{r-2}|_{S^{r-2}}\Big)^{-1}(\lambda).$  By induction hypothesis applied on $S^{r-2}$, we get that for $\mu \in \Big(\Phi_{2}^{r-2}|_{S^{r-2}}\Big)^{-1}(\lambda),$ $\deg P_\mu <r-2.$ Thus, we have 
$$\deg \sum_{\mu \in \big(\Phi_2^{r-2}|_{S^{r-2}}\big)^{-1}(\lambda)} (2m - 2 - |\mu|) P_\mu(m-1)  < (r-2 ) + 1 = r-1.$$

On the other hand, we have \[
\deg \bigl[x_0^{\,2m-|\lambda|}x^\lambda\bigr]
\bigl(x_0x_1\,\delta(\Gamma_{2m-2, 2m-r})\bigr)
= \]
\[\deg \sum_{\nu \in (\Phi_1^{\,r-1}|_{U^{\,r-1}_{\mathrm{no1}}})^{-1}(\lambda)}
\;
\sum_{(\mu,u)\in (\Psi_1^{\,r-2})^{-1}(\nu)}
\varepsilon(\mu,u)\,P_\mu(m-1) < r-1,\]
because the $\mu$'s appearing inside belong to $ S^{r-2}$ and by induction hypothesis $\deg P_{\mu} \leq r-2 < r-1.$ Moreover, we have
\begin{align*}
&\deg \big[x_0^{2m-|\lambda|}x^\lambda\big]\Big( x_0x_1\Gamma_{2m-2, 2m-1-r}\Big)\\ =& \deg [x_0^{2m-|\lambda|}x^\lambda] \Big(\sum_{\mu \in S^{r-1}} P_\mu(m-1) x_0^{2m-1-|\mu|}x_1x^{\mu}\Big)\\ =& \deg \sum_{\mu \in \big(\Phi_1^{r-1}|_{\{\mu \in S^{r-1} : 1 \notin \mu\}}\big)^{-1}(\lambda)} P_\mu(m-1) \leq  \max_{\mu \in \big(\Phi_1^{r-1}|_{\{\mu \in S^{r-1} : 1 \notin \mu\}}\big)^{-1}(\lambda)}(\deg P_\mu(m-1)).
\end{align*}
If $r$ is odd, then for $ \mu \in \big(\Phi_1^{r-1}|_{\{\mu \in S^{r-1} : 1 \notin \mu\}}\big)^{-1}(\lambda), \; \deg P_\mu(m) < r-1$ unless $\lambda_{(r-1) \text{-packed}} \in \big(\Phi_1^{r-1}|_{\{\mu \in S^{r-1},~ 1 \notin \mu\}}\big)^{-1}(\lambda).$ The latter happens precisely when $\lambda$ itself is packed from \eqref{eq:prop2a}, which is a contradiction. If $r$ is even, it follows from \eqref{eq:prop2b} that $P_{\lambda_{(r-1)\text{-packed}}}$ does not appear in the sum. Thus, the sum has degree $ < r-1.$ 

Therefore, we have just proved that right-hand side of \eqref{eq:RHS} has degree $<r-1$ for not packed $\lambda.$ A finite integration leads to $\deg P_\lambda(m) < r,$ for not packed  $\lambda$. Thus, if $\deg P_\lambda(m) = r$, we must have that $\lambda$ is packed. 
\end{proof}

It might technically happen  that the degree $r$ is simply unattainable. We will demonstrate that packed shapes have degree $r$ by inductively computing the coefficient of $m^r$ in $P_{\lambda_{r-\text{packed}}}$.

\sssbegin{Lemma}\label{lem:A17}
    Let $\lambda$ be a $r$-packed shape and let $f(r)$ denote the leading coefficient of $P_\lambda(m)$. Then, we have $f(r)=\frac{1}{\lfloor r/2\rfloor!}.$
\end{Lemma}
\begin{proof}
    We will show the Lemma by induction on $r$. We have $f(1) = f(2) = 1.$
The degree of
$$\sum_{\nu \in (\Phi_1^{\,r-1}|_{U^{\,r-1}_{\mathrm{no1}}})^{-1}(\lambda)}
\;
\sum_{(\mu,u)\in (\Psi_1^{\,r-2})^{-1}(\nu)}
\varepsilon(\mu,u)\,P_\mu(m-1) $$
is at maximum $r-2,$ so this term does not contribute to the coefficient of $m^{r-1}$.

\underline{The case where $r = 2l$ is even}. Then $r-1$ is odd and $r-2$ is even. In this case, the packed $r-$shape is $\lambda = (2^{\frac{r}{2}})$. The leading coefficient of 

$$\sum_{\mu \in \big(\Phi_2^{r-2}|_{S^{r-2}}\big)^{-1}(\lambda)} (2m - 2 - |\mu|) P_\mu(m-1) + \sum_{(\mu, u) \in \big(\Psi_2^{r-2}|_{S^{r-2}}\big)(\lambda)} P_\mu(m-1) $$
comes out from $(2m-2-|\mu|)P_\mu(m-1)$ where $\mu$ is $(r-2)$-packed shape. This term has a leading coefficient $2f(r-2).$

Notice that the term
$$\sum_{\mu \in \big(\Phi_1^{r-1}|_{\{\mu \in S^{r-1} : 1 \notin \mu\}}\big)^{-1}(\lambda)} P_\mu(m-1)$$
contributes precisely only when $\mu$ is $(r-1)$-packed. But in that case, such $\mu$ must contain $1$, so it does not even appear in the summand. Thus, we would have that the leading term in the right-hand side of \eqref{eq:RHS} is simply $2f(r-2)$ in that case. Since we have a finite difference applied to $P_\lambda(m)$, this means that
$$r \cdot f(r) = 2f(r-2).$$
It follows that $f(r) = \frac{2}{r}f(r-2)$
when $r$ is even. From the induction hypothesis, $$f(r-2) = \frac{1}{\lfloor \frac{r-2}{2} \rfloor!} = \frac{\lfloor \frac{r}{2} \rfloor}{\lfloor \frac{r}{2}\rfloor!} = \frac{r}{2(r/2)!}.$$
Thus, 
$$f(r) = \frac{2}{r} \cdot \frac{r}{2} \cdot \frac{1}{\lfloor\frac{r}{2}\rfloor!} = \frac{1}{\lfloor \frac{r}{2}\rfloor!},$$
finishing the induction step in this case.

\underline{The case where $r = 2l+1$ is odd}. As in the previous case, the leading coefficient of 
$$\sum_{\mu \in \big(\Phi_2^{r-2}|_{S^{r-2}}\big)^{-1}(\lambda)} (2m - 2 - |\mu|) P_\mu(m-1) + \sum_{(\mu, u) \in \big(\Psi_2^{r-2}|_{S^{r-2}}\big)(\lambda)} P_\mu(m-1) $$
comes out from $(2m-2-|\mu|)P_\mu(m-1)$ where $\mu$ is $(r-2)$-packed shape. This term has a leading coefficient of $2f(r-2).$ Unlike the case where $r$ is even, the term $$\sum_{\mu \in \big(\Phi_1^{r-1}|_{\{\mu \in S^{r-1} : 1 \notin \mu\}}\big)^{-1}(\lambda)} P_\mu(m-1)$$ can contribute. The contribution is precisely the leading coefficient of $P_\mu(m-1)$ where $\mu = (2^{\frac{r-1}{2}})$, the $(r-1)-$packed shape which does not contain $1.$ The contribution here is $f(r-1).$ Thus, we obtain the following relation for odd $r:$
$$rf(r) = 2f(r-2) + f(r-1).$$
Since $f(r-2) = \frac{1}{\lfloor \frac{r-2}{2}\rfloor!} =  \frac{1}{(\frac{r-3}{2})!}$ and $f(r-1) = \frac{1}{(\frac{r - 1}{2})!}$, from the induction hypothesis and using that $r$ is odd, we have
\begin{align*}
f(r) &= \frac{2}{r} \cdot \frac{1}{(\frac{r-3}{2})!} + \frac{1}{r \cdot (\frac{r-1}{2})!} = \frac{1}{(\frac{r-3}{2})!} \Big(\frac{2}{r} + \frac{1}{r\cdot \frac{r-1}{2}}\Big)= \frac{1}{(\frac{r-1}{2})!},
\end{align*}
which finishes the induction proof.
\end{proof}

\sssbegin{Remark} Lemmas \ref{lem:A16} and \ref{lem:A17} together imply that $\deg P_\lambda(m) = r$ if and only if $\lambda$ is packed. 
\end{Remark}

\subsection{The proof of the Proposition \ref{prop:p-divides-Gamma2pp}}

    We are finally able to show that $\Gamma_{2p,p}\equiv 0 \mod p$. The polynomials $P_\lambda(m)$ are integer valued. Moreover, $P_\lambda(0) = 0$ for any $\lambda \in S^r$ with $r \geq 1$ and $\deg P_\lambda(m) \leq r$ for $\lambda \in S^r.$

Set $r = p$. For $\lambda \in S^p$, the polynomials $P_\lambda(m)$ admit a binomial coefficient basis, i.e. there exists $c_0, c_1, \hdots, c_{\deg P_\lambda(m)} \in \mathbb{Z}$ such that
$$P_\lambda(m) = \sum_{i=0}^{\deg P_\lambda(m)} c_i \cdot \binom{m}{i}.$$
Suppose that $\lambda$ is not a packed shape. Then, $\deg P_\lambda(m) < p.$ For $m = 0$, we have
$$0 = P_\lambda(0) = \sum_{i=0}^{\deg P_\lambda(m)} c_i \binom{m}{i} = c_0.$$
Thus, $c_0 = 0.$ Let $m = p$ and consider the expression modulo $p$. For any $1 \leq i \leq p - 1$, we have $\binom{p}{i} \equiv 0$ and $\deg P_\lambda < p.$ It follows that
$$P_\lambda(p) \equiv \sum_{i=0}^{\deg P_\lambda(p)} c_i \binom{p}{i} \equiv 0 \pmod{p}.$$
On the other hand, in the case where $\lambda$ is packed, we have $\deg P_\lambda(p) = p.$ In particular, we have
$$P_\lambda(p) \equiv \sum_{i=0}^p c_i \binom{p}{i} \equiv c_p \mod p.$$
However, the leading coefficient $l$ of $P_\lambda(p)$ and $c_{\deg P_\lambda(p)}$ satisfy the relation $l = \frac{c_{\deg P_\lambda(p)}}{(\deg P_\lambda(p))!}.$ Since we have $l = \frac{1}{\lfloor \frac{p}{2}\rfloor!} = \frac{1}{(\frac{p-1}{2})!}$, it follows that $c_p = \frac{p!}{(\frac{(p-1)}{2})!} \equiv 0 \pmod p.$  Thus, we have  $P_\lambda(p) \equiv c_p \equiv 0 \mod p$ for any packed $\lambda$.

Finally, $\Gamma_{2p,p} = \displaystyle\sum_{\lambda \in S^p} P_\lambda(p) x_0^{2p - |\lambda|}x^\lambda \equiv 0 \pmod p$ and Proposition \ref{prop:p-divides-Gamma2pp} is proved.\quad\quad\quad\quad$\square$


\end{document}